\documentclass[11pt]{article}

\usepackage{comment}
\usepackage{amsmath, amsthm, amsfonts}
\usepackage{listings}
\usepackage{amsfonts}
\usepackage{amssymb}
\usepackage{ulem}
\usepackage{float}
\usepackage{anysize}
\usepackage{graphicx}
\usepackage{mathrsfs} 
\usepackage{hyperref}
\usepackage{cite}
\usepackage{color}
\usepackage{subcaption}
\usepackage[round, authoryear]{natbib} 
\usepackage{color}
\usepackage{cancel}
\usepackage{bbm}
\usepackage{soul}
\usepackage{ulem}
\usepackage{cancel}
\usepackage{tikz}
\usepackage{lscape}
\usepackage{lineno} 

\usepackage{fullpage} 
\usepackage{multirow}
\usepackage{lipsum}
\usepackage{mwe}
\usepackage{xcolor} 
\usepackage{colortbl}
\usepackage{lscape}

\newtheorem{thm}{Theorem}[section]

\newtheorem{lem}[thm]{Lemma}
\newtheorem{prop}[thm]{Proposition}
\theoremstyle{definition}

\theoremstyle{remark}
\newtheorem{rem}{Remark}
\theoremstyle{remark}
\newtheorem{Assum}{Assumption}

\usepackage{authblk}
\title{Asymptotic Consistency of a Weighted Nadaraya-Watson Estimator Driven by Fractional Brownian Motion at Random Times.}
\author[1,2]{Tania Roa \thanks{tania.roa@uai.cl}}
\author[1]{Meryem Slaoui \thanks{meryem.slaoui@sgcib.com}}
\author[3]{Soledad Torres \thanks{soledad.torres@uv.cl}}

\affil[1]{Universidad Adolfo Ibáñez, Faculty of Engineering \& Sciences. Viña del Mar, Chile}
\affil[2]{Data Observatory Foundation, ANID Technology Center No. DO210001, Eliodoro Yáñez 2990, 7510277, Providencia, Santiago, Chile}
\affil[3]{CIMFAV - Facultad de Ingeniería. Universidad de Valparaíso, Valparaíso, Chile.}
\usepackage{microtype} 

\begin{document}
\maketitle

\begin{abstract}
Nonparametric regression conventionally assumes equispaced observations, a premise often violated by empirical data with irregular sampling and long-memory noise. This study establishes the asymptotic mean square ($L^2$) consistency of a weighted Nadaraya-Watson (wN-W) estimator for stochastic processes sampled at random times and perturbed by fractional Brownian motion (fBm). Leveraging the occupation time formula and Fernique's Theorem, we derive exact parametric restrictions on bandwidth and sampling rates ensuring convergence. We validate these theoretical findings via extensive simulations evaluating the mean squared error under varying noise persistence. Crucially, we explicitly benchmark the performance of the estimator under random sampling times against classical deterministic grids, while assessing the relative efficiency of quartic versus Gaussian kernels. Finally, we introduce an interactive Shiny application to dynamically visualize the admissible parameter space, and we empirically validate our theoretical framework by applying the estimator to the Vostok ice core paleoclimatological dataset. This real-world application demonstrates the practical capacity of the methodology to extract macroscopic trends while robustly avoiding severe overfitting under random sampling times.\\

\vspace{9pt}
\end{abstract}

\section{Introduction}

Nonparametric regression provides a flexible framework to model complex dependence structures without imposing rigid functional forms, with widespread applications ranging from time series \cite{tjostheim1994} and survey analysis \cite{ling2018} to image processing \cite{andreetto2007} ; \cite{hu2010}). Among these techniques, the Nadaraya-Watson (N-W) estimator \cite{nadaraya1964} and \cite{watson1964} remains a cornerstone. While the N-W estimator has been successfully adapted for non-stationary time series (\cite{soukarieh2023}), fuzzy data \cite{hesamian2023}, and missing data (\cite{cheung2023}), its asymptotic performance heavily depends on the choice of the kernel (\cite{wasserman2006}) and the bandwidth $h \equiv N^{-\alpha}$ ($\alpha > 0$). In this paper, we focus on establishing the exact parametric conditions for $\alpha$ to ensure the estimator's convergence under complex temporal and noise structures.
A significant challenge in nonparametric estimation arises when observations exhibit long-memory dependence. While discrete-time long memory is often modeled via ARFIMA processes (\cite{granger1980}; \cite{guegan2005}), continuous-time phenomena with long-range dependence are best captured by fractional Brownian motion (fBm). The fBm is a centered Gaussian process with stationary increments and a covariance function defined by $\mathbb{E}[B_t^H B_s^H] = \frac{\sigma^2}{2}[t^{2H} + s^{2H} - \vert{}t-s\vert{}^{2H}]$. Extensive research has addressed parametric estimation in regression models driven by fBm or mixed fBm (\cite{wang2015}; \cite{chen2021}). In the nonparametric realm, earlier works have explored the impact of long memory \cite{estevez2003} and applied nonparametric methods to estimate trends in stochastic differential equations driven by mixed fBm (\cite{prakasa2019}).

Moving beyond the traditional constraints of equidistant observations, such as those assumed in the nonparametric long-memory models of \cite{sued2013}, we establish a generalized Nadaraya-Watson estimation strategy tailored for randomly sampled data. Irregular sampling has been addressed in spectral analysis (\cite{stoica2011}) and spatial statistics (\cite{lu2014}), and recently, nonparametric drift estimation for SDEs has been studied using multiple i.i.d. paths (\cite{han2023}). However, our approach focuses on sampling times governed by a specific probability distribution. Specifically, we follow the renewal process framework recently developed by \cite{araya2023} for parametric least squares (LS) estimators, adapting it to establish the consistency of the N-W estimator when the noise is driven by fBm. \\

While adapting these theoretical frameworks establishes the necessary asymptotic bounds for the bandwidth parameter $\alpha$, the transition to highly irregular, non-equidistant random sampling times introduces severe empirical vulnerabilities. If the optimal theoretical constraints for $\alpha$ are ignored or violated in practical applications, the N-W estimator risks fitting the high-frequency stochastic noise of the fBm rather than the underlying structural trend. To demonstrate that our derived regularity conditions are not merely incremental mathematical constraints, but rather indispensable operational boundaries against severe overfitting, this paper also provides a comprehensive empirical validation using the Vostok ice core paleoclimatological dataset. By bridging the theoretical limits of $\alpha$ with a complex real-world irregular sampling scenario, we illustrate the direct practical implications of our consistency theorems.


\noindent The remainder of this paper is organized as follows. Section \ref{prelim} establishes the mathematical framework, detailing the regression model, the renewal process governing the irregular sampling times, and the long-memory noise structure. It also outlines the necessary structural assumptions—particularly regarding the kernel functions and the fractional moments of the temporal gaps—and formally defines the weighted Nadaraya-Watson (wN-W) estimator. Section \ref{main-res} presents the core theoretical contributions, introducing the sequence of lemmas and the main theorem that establish the exact parametric conditions for the $L^{2}$-consistency of the estimator. Section \ref{proofs} is devoted to the rigorous proofs of these findings. Here, we deploy advanced tools from stochastic analysis, such as the occupation time formula and Fernique's Theorem, to bound the discretization errors and control the asymptotic behavior of the fractional noise. Section \ref{simulations} provides a comprehensive simulation study to validate our theoretical claims. We analyze the finite-sample performance and the Mean Squared Error (MSE) decay across different bandwidth regimes ($\alpha$) and Hurst exponents ($H$), and we introduce an interactive open-source web application designed to dynamically visualize the admissible parametric space. Finally, Section \ref{real-data} presents a comprehensive empirical application using the Vostok ice core paleoclimatological dataset. This concluding section demonstrates the practical applicability of our estimator on real-world non-equidistant data and highlights the severe overfitting that occurs when the derived topological constraints are violated.


\section{Preliminaries}
\label{prelim}

\subsection{The Model and Random Times}
We consider the following nonparametric regression model sampled at random times, inspired by \cite{sued2013}:
\begin{equation}\label{no-p-rt}
Y_{\tau_{i}}=r(x_{\tau_{i}})+\epsilon_{\tau_i}^{H_{2}}, \quad i=0,1,...,N(1);
\end{equation}
with $\tau_{0}=0$. The regression function $r(\cdot)$ is deterministic and unknown. The covariates are given by $x_{\tau_{i}}=B_{\tau_{i}}^{H_{1}}$, where $B^{H_{1}}$ is a fractional Brownian motion (fBm) with Hurst index $1/2 < H_{1} < 1$, assumed independent of the noise $\epsilon_{\tau_i}^{H_{2}}$. 

The observation interval is set to $[0,1]$, containing a random number of observations $N(1)$. The sequence of random sampling times $\tau=\{\tau_{i}\}_{i=0}^{N(1)}$ is strictly increasing, defined by the renewal process $\tau_{i}=\sum_{j=0}^{i}t_{j}$ for $i\ge 0$. Here, $\{t_{j}\}_{j\ge 0}$ is an i.i.d. sequence of inter-sampling times with probability density $g_{t_{i}}$ (\cite{araya2023}). Consequently, the cumulative time $\tau_{i}$ admits a density $G_{\tau_{i}}$. Generic constants depending on the moments of $t_i$ will be denoted by $C$ or $C_\kappa$ throughout the paper.

For the estimation procedure, we employ a kernel function $K(\cdot)$. We consider both the standard Gaussian kernel and bounded support kernels (e.g., triangle, Epanechnikov, and quartic kernels; see \cite{hardle1990}). For bounded kernels, we assume $K(x) = g_k(x)\mathbb{I}_{[-1,1]}(x) \le C_K$, where $C_K$ is a known bounding constant.

\subsection{The Noise Structure}
The model noise $\epsilon_{\tau_{i}}^{H_{2}}$ is driven by the increments of an independent fBm $B^{H_{2}}$ with $1/2 < H_{2} < 1$. To ensure that the error exhibits asymptotically constant variance, we must normalize these fractional increments. By conditioning on $\tau_{i}=z_{i}$ and $t_{i+1}=w_{1}$, the variance of the increment is bounded as follows:
\begin{align}
&Var(B_{\tau_{i+1}}^{H_{2}}-B_{\tau_{i}}^{H_{2}}) = \mathbb{E}[\mathbb{E}[(B_{z_{i}+w_{1}}^{H_{2}}-B_{z_{i}}^{H_{2}})^{2}]|\tau_{i}=z_{i},t_{i+1}=w_{1}] \nonumber \\
&= \int_{0}^{\infty}\int_{0}^{\infty}w_{1}^{2H_{2}}G_{\tau_{i}}(z_{i})g_{t_{i+1}}(w_{1})dz_{i}dw_{1} 
= \mathbb{E}[t_{i+1}^{2H_{2}}] \le (\mathbb{E}[t_{i+1}^{2}])^{H_{2}} \le \frac{\kappa_{1}^{H_{2}}}{N^{aH_{2}}}. \label{constante-m}
\end{align}
Based on this bound, we define the normalized noise for our model as $\epsilon_{\tau_{i}}^{H_{2}}=\frac{N^{\frac{aH_{2}}{2}}}{\kappa_{1}^{\frac{H_{2}}{2}}}(B_{\tau_{i+1}}^{H_{2}}-B_{\tau_{i}}^{H_{2}})$.

\subsection{Assumptions}  \label{asump}
In this section, we outline the structural assumptions that underpin our asymptotic analysis.

\begin{Assum}[Kernel and Function Regularity]\label{assump:regularity}
(i) The kernel $K(\cdot)$ is a bounded ($\sup_x K(x) \le C_K < \infty$), square-integrable function satisfying $\int_{\mathbb{R}} K(y)dy = 1$. (ii) The regression function $r(\cdot)$ is square-integrable ($\int |r(x)|^2 dx < \infty$) and locally $\gamma_r$-H\"older continuous, with $0 < \gamma_r \le 1$. (iii) The fractional Brownian motions exhibit H\"older continuity of order $\delta \in (0, H_1)$ (\cite{nourdin2012}). Similar H\"older conditions ($\gamma_f$) apply to any test function $f$ used throughout the proofs.
\end{Assum}

\begin{Assum}[Fractional Inverse Moments]\label{assump:inverse}
For a strictly positive random variable $Y$ and constant $\mu>0$, the fractional inverse moment is finite, i.e., $\mathbb{E}[Y^{-\mu}] < C_\mu < \infty$. Following \cite{cressie1986}, a sufficient condition is $\int_{0}^{\infty}t^{\mu-1}M_{Y}(-t)dt<\infty$, where $M_Y$ denotes the moment generating function.
\end{Assum}

\begin{rem}
Assumption \ref{assump:inverse} holds for typical inter-arrival time distributions. For instance, if $Y_i$ are exponential random variables with parameter $\lambda=N$, the convolution $\tilde{Y}_{k}=\sum_{i=1}^{k}Y_{i}$ yields $M_{\tilde{Y_{k}}}(t)=(N/(N-t))^{k}$, which ensures a finite inverse fractional moment via hypergeometric functions. Other suitable distributions include the Rayleigh distribution, standard in wind speed analysis (\cite{david2014}), and the Gamma/Gompertz distribution, prevalent in mortality risk models (\cite{missov2013}).
\end{rem}

\begin{Assum}[Renewal Process Dynamics]\label{assump:renewal}
Following Araya et al. (2023), the i.i.d. inter-sampling times $t_i$ satisfy the following moment conditions:
\begin{equation} \label{hipotesis3}
\mathbb{E}[t_{i}]=\frac{1}{N}, \quad \mathbb{E}[t_{i}^{2}]=\frac{\kappa_{1}}{N^{a}}, \quad \text{and} \quad \mathbb{E}[t_{i}^{4}]=\frac{\kappa_{2}}{N^{b}},
\end{equation}
where $a \le 2$ and $0 < b < 2a$ (by the Cauchy-Schwarz inequality). Under these conditions, the counting process $N(1)$ and the temporal gaps exhibit the following asymptotic properties as $N \rightarrow \infty$:
\begin{enumerate}
    \item[(R1)] $\tau_{N-1} \rightarrow 1$ and $\frac{N(1)}{N} \rightarrow 1$. Both convergences hold in probability for $a>1, b>1$, and almost surely for $a>3/2, b>2$.
    \item[(R2)] For any $\epsilon > 0$, $\tau_{N+N\epsilon}-\tau_{N} \rightarrow \epsilon$ almost surely, provided $a>3/2$ and $b>2$.
\end{enumerate}
\end{Assum}

Finally, to control the rate of convergence of the $N$-th observation time to 1, we rely on the following bounds derived from Assumption \ref{assump:renewal}:
\begin{equation}\label{fourth}
\mathbb{E}[|\tau_{N}-1|^{2}]\le\frac{C}{N^{a-1}} \quad \text{and} \quad \mathbb{E}[|\tau_{N}-1|^{4}]\le\frac{C}{N^{(2a-2)\wedge(b-1)}}.
\end{equation}

\section{Main Results} \label{main-res}

To establish the consistency of the wN-W estimator, we evaluate its asymptotic behavior under both the random sample size $N(1)$ and a fixed sample size $N$. Let $\Delta\tau_{i}=\tau_{i+1}-\tau_{i}=t_{i+1}$. Following \cite{cai2001}, we define the auxiliary functions for the numerator and denominator for any integer $n \ge 1$:
\begin{equation} \label{eq:AnDn}
A_{n}(x)=\sum_{i=0}^{n-1}\Delta\tau_{i}Y_{\tau_{i}}K((x-B_{\tau_{i}}^{H_{1}})N^{\alpha}), \quad D_{n}(x)=\sum_{i=0}^{n-1}\Delta\tau_{i}K((x-B_{\tau_{i}}^{H_{1}})N^{\alpha}).
\end{equation}
The wN-W estimator is thus expressed as $\hat{r}_n(x) = A_n(x)/D_n(x)$. The total estimation error can be bounded by decomposing it into the discrepancy caused by the random stopping time and the fixed-sample estimation error: $|\hat{r}_{N(1)}(x)-r(x)| \le |\hat{r}_{N(1)}(x)-\hat{r}_{N}(x)| + |\hat{r}_{N}(x)-r(x)|$. 

Furthermore, the fixed-$N$ estimator is analytically split into a noise-free approximation term $T_{1}^{(N)}(x)$ and a fractional noise perturbation term $T_{2}^{(N)}(x)$, such that $\hat{r}_N(x) = T_{1}^{(N)}(x) + T_{2}^{(N)}(x)$.

\begin{thm} \label{thm:main}
Let the structural conditions in Assumptions \ref{assump:regularity}--\ref{assump:renewal} hold with $1 < a < 2$ and $2 < b < 2a$. If the bandwidth parameter $\alpha$ satisfies
\begin{equation} \label{eq:alpha_bounds}
\begin{split}
\alpha < \min\Bigg( &\frac{a(H_1\gamma_r+1)-2}{2(\gamma_r+1)}, \frac{b(H_1\gamma_r+1)-2}{4(\gamma_r+1)}, \frac{a(H_1\gamma_K+1)-2}{2(\gamma_K+1)}, \frac{b(H_1\gamma_K+1)-2}{4(\gamma_K+1)}, \\ 
&\frac{a-1}{2(1 + \gamma_K)}, \frac{b-1}{4(1 + \gamma_K)}, \alpha_K^* \Bigg),
\end{split}
\end{equation}
where $\alpha_K^* = \min\left( \frac{2}{3}\left(\frac{b}{2}-1-H_2\left(a-\frac{b}{2}\right)\right), \frac{2}{3}\left(2-\frac{H_1}{12}-aH_2\right) \right)$ for bounded support kernels, and $\alpha_K^* = a\left(2 - H_2 - \frac{5H_1}{48}\right) - 2$ for the Gaussian kernel, then the wN-W estimator is $L^2$-consistent:
$$ \hat{r}_{N(1)}(x) \xrightarrow{L^2} r(x) \quad \text{as } N \to \infty. $$
\end{thm}

The proof of Theorem \ref{thm:main} relies on bounding the two terms of the error decomposition via the following propositions.

\begin{prop} \label{prop:fixedN}
Under the conditions of Theorem \ref{thm:main}, the fixed-sample estimator is $L^{2}$-consistent, i.e., $\hat{r}_{N}(x) \xrightarrow{L^2} r(x)$ as $N \rightarrow \infty$.
\end{prop}

\begin{prop} \label{prop:randomN}
Under the conditions of Theorem \ref{thm:main}, the discrepancy between the estimators evaluated at $N(1)$ and $N$ vanishes in mean square: $|\hat{r}_{N(1)}(x)-\hat{r}_{N}(x)| \xrightarrow{L^2} 0$ as $N \rightarrow \infty$.
\end{prop}

Proposition \ref{prop:fixedN} is proven by controlling the convergence of $T_{1}^{(N)}(x)$ and $T_{2}^{(N)}(x)$ through a sequence of auxiliary lemmas. Let the step function be $c(s)=\sum_{i\ge 0}\tau_{i}\mathbb{I}_{[\tau_{i},\tau_{i+1})}(s)$.

\begin{lem} \label{lem:step}
For any $q \in (0,1)$, the $L^2$ norm of the temporal discretization is bounded by $\mathbb{E}[(\int_{0}^{1}|c(s)-s|^{q}ds)^{2}]\le C N^{-[(\frac{b(q+1)}{2}-1)\wedge(a(q+1)-2)]}$.
\end{lem}

\begin{lem} \label{lem:integrals}
Let $f$ satisfy the H\"older conditions in Assumption \ref{assump:regularity}. If $\alpha < \min(\frac{a(H_{1}\gamma_{f}+1)-2}{2(\gamma_{f}+1)}, \frac{a-1}{2})$, the discrete sum approximating the denominator converges to the continuous integral:
$$ R_{N} := N^{\alpha}\left\{\sum_{i=0}^{N-1}\Delta\tau_{i}f((x-B_{\tau_{i}}^{H_{1}})N^{\alpha})-\int_{0}^{1}f((x-B_{s}^{H_{1}})N^{\alpha})ds\right\} \xrightarrow{L^2} 0. $$
Analogously, under the full constraints of Theorem \ref{thm:main}, the numerator approximation error $L_N$ (incorporating the functions $r(\cdot)$ and $K(\cdot)$) also converges to $0$ in $L^2$.
\end{lem}

\begin{lem}\label{lem:localtime}
Assume $f$ satisfies $\int_{\mathbb{R}}|f(z)||z|^{\delta}dz<\infty$ for $\delta \in (0, \frac{1-H_1}{2H_1})$. As $N \to \infty$, the fractional integral converges to the local time:
$$ N^{\alpha}\int_{0}^{1}f((x-B_{c(s)}^{H_{1}})N^{\alpha})ds \xrightarrow{L^2} d_{1}L^{H_{1}}(1,x), $$
where $d_{1}=\int_{\mathbb{R}}f(y)dy$ and $L^{H_{1}}(1,x)$ is the local time of the fBm at time $1$ and spatial level $x$ (\cite{sued2013}). This result naturally extends to the product $r(B_{c(s)}^{H_1})f(\cdot)$, converging to $d_1 r(x) L^{H_1}(1,x)$.
\end{lem}

\begin{rem}
The strict positivity of the local time ($L^{H_{1}}(1,x)>0$ a.s.) for any spatial level $x$ in the interior of the process range ensures that the denominator of the wN-W estimator is asymptotically well-defined.
\end{rem}

\begin{lem} \label{lem:noise}
Let $K$ be a valid kernel (Assumption \ref{assump:regularity}). If $\alpha$ satisfies the kernel-specific upper bounds defined by $\alpha_K^*$ in Theorem \ref{thm:main}, the noise component $Q_N = N^{\alpha}\sum_{i=0}^{N-1}\Delta\tau_{i}\epsilon_{\tau_{i}}^{H_{2}}K((x-B_{\tau_{i}}^{H_{1}})N^{\alpha})$ converges to $0$ in $L^2$.
\end{lem}



\section{Proofs of Propositions, Lemmas and main Theorem} \label{proofs}

This section presents the proofs of Lemmas \ref{lem:step}, \ref{lem:integrals}, \ref{lem:localtime}, and \ref{lem:noise}, which together form the basis of the proof of Proposition \ref{prop:fixedN}. Thereafter, we proceed to the proof of Proposition \ref{prop:randomN}, which, when joined with Proposition \ref{prop:fixedN}, yields the proof of the main result, Theorem \ref{thm:main}.

\subsection*{Proof of Lemma \ref{lem:step}}
\begin{proof}
Let us consider the following figure illustrating $c(s)$, $\tau_{i}$, $\tau_{i+1}$, and $s$ along the line.

\begin{figure}[h!]
\centering
\includegraphics[scale=0.25]{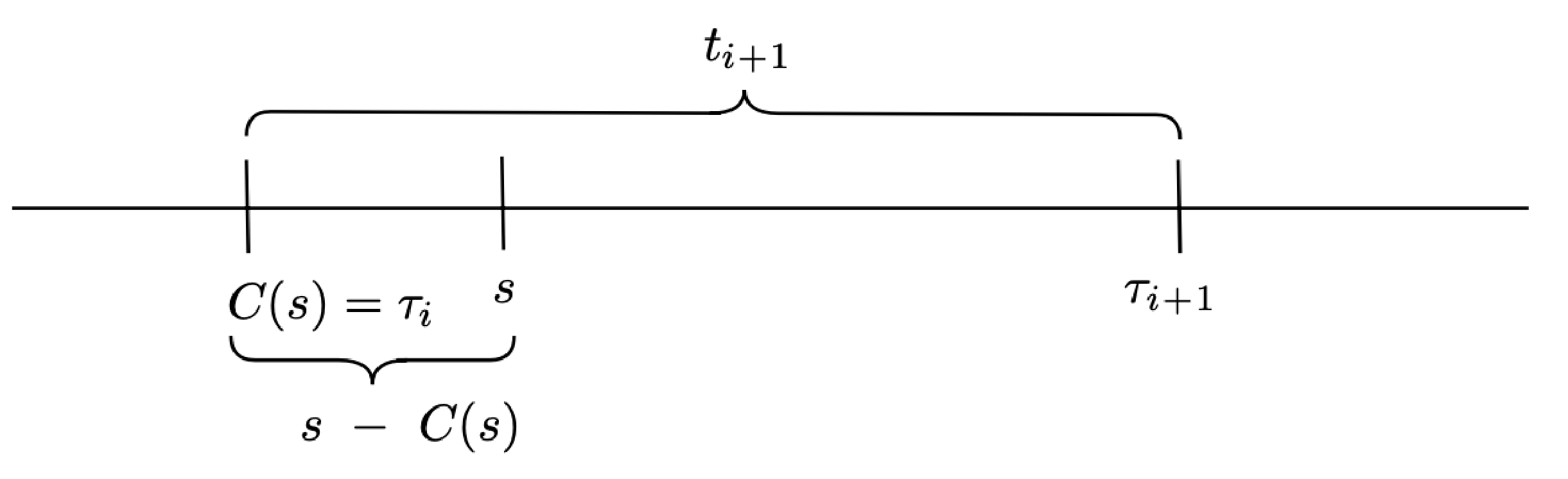}
\caption{Representation of $c(s)$}
\end{figure}

For $0 < q < 1$, by partitioning the integral and applying the inequality $(x+y)^2 \le 2x^2 + 2y^2$, we have:
\begin{align*}
&\mathbb{E} \left[ \left( \int_{0}^{1} \left| c(s) - s \right|^{q} ds \right)^{2} \right] = \mathbb{E} \left[ \left( \sum_{i=0}^{N-1} \int_{\tau_{i}}^{\tau_{i+1}} \left| \tau_{i} - s \right|^{q} ds + \int_{\tau_{N}}^{1} \left| \tau_{N} - s \right|^{q} ds \right)^{2} \right] \\
&\leq 2 \mathbb{E} \left[ \left( \sum_{i=0}^{N-1}  \left| t_{i+1} \right|^{q+1} \right)^{2} \right] + 2 \mathbb{E} \left[  \left| 1 - \tau_{N}  \right|^{2q+2}  \right] \\
&= 2 \sum_{i=0}^{N-1} \mathbb{E} \left[ \left| t_{i+1} \right|^{2q+2} \right] + 4 \sum_{0 \leq i < j \leq N-1} \mathbb{E} \left[ \left| t_{i+1} \right|^{q+1} \right]  \mathbb{E} \left[ \left| t_{j+1} \right|^{q+1} \right] + 2 \mathbb{E} \left[  \left| 1 - \tau_{N}  \right|^{2q+2}  \right].
\end{align*}

From the bounds on the final time \eqref{fourth} and applying Jensen's inequality to the temporal moments, we obtain:
\begin{align*}
&\leq 2 \sum_{i=0}^{N-1} \left( \mathbb{E} \left[ \left| t_{i+1} \right|^{4} \right]\right)^{\frac{q+1}{2}} + 4 \sum_{0 \leq i < j \leq N-1} \left( \mathbb{E} \left[ \left| t_{i+1} \right|^{2} \right] \right)^{\frac{q+1}{2}}  \left( \mathbb{E} \left[ \left| t_{j+1} \right|^{2} \right]  \right)^{\frac{q+1}{2}} + 2 \left( \mathbb{E} \left[  \left| 1 - \tau_{N}  \right|^{4} \right] \right) ^{\frac{q+1}{2}} \\
&= 2 \sum_{i=0}^{N-1} \left( \dfrac{\kappa_{2}}{N^{b}} \right)^{\frac{q+1}{2}} + 4 \sum_{0 \leq i < j \leq N-1}  \left( \dfrac{\kappa_{1}}{N^{a}} \right)^{q+1} + 2 \left( \mathbb{E} \left[  \left| 1 - \tau_{N}  \right|^{4} \right] \right) ^{\frac{q+1}{2}} \\
&\leq \dfrac{C}{N^{\frac{b(q+1)}{2} - 1}} +  \dfrac{C}{N^{a(q+1) - 2}} + \dfrac{C}{N^{(b-1) \left(\frac{q+1}{2} \right) \wedge (2a-2)  \left(\frac{q+1}{2} \right)} } \\
&\leq \dfrac{C}{N^{\left[ \frac{b(q+1)}{2} - 1 \right] \wedge \left( a(q+1) - 2 \right) }}.
\end{align*}
\end{proof}

\subsection*{Proof of Lemma \ref{lem:integrals}}
\begin{proof}
As established in the Preliminaries section, an important result is the convergence of the $N^{th}$ random time, $\tau_{N}$ to 1 if $a > 3/2$ and $b>2$. First, we prove the convergence of the denominator sequence $R_N$. Rewriting its Riemann-Stieltjes sum:
\begin{align}
\sum_{i=0}^{N-1} \left( \Delta \tau_{i} \right)  f \left( \left( x - B_{\tau_{i}}^{H_1} \right) N^{\alpha} \right) &= \int_{0}^{\tau_{N}} f \left( \left( x - B_{c(s)}^{H_1} \right) N^{\alpha} \right) ds \nonumber \\
&= \int_{0}^{1} f \left( \left( x - B_{c(s)}^{H_1} \right) N^{\alpha} \right) ds - \int_{\tau_{N}}^{1} f \left( \left( x - B_{c(s)}^{H_1} \right) N^{\alpha} \right) ds, \label{sum-int}
\end{align}
where $c(s) = \tau_{i} \mathbbm{1}_{[\tau_{i}, \tau_{i+1})}(s)$. Therefore, $R_{N}$ can be decomposed into an internal discretization term and a tail error:
\begin{equation} \label{RN}
R_{N} = N^{\alpha} \left\lbrace \int_{0}^{1} \left[ f \left( \left( x - B_{c(s)}^{H_1} \right) N^{\alpha} \right) - f \left( \left( x - B_{s}^{H_1} \right) N^{\alpha} \right) \right] ds +  \int_{\tau_{N}}^{1} f \left( \left( x - B_{c(s)}^{H_1} \right) N^{\alpha} \right) ds  \right\rbrace. 
\end{equation}

Given the H\"older continuity and boundedness of $f(\cdot)$, for $N$ large enough, we bound the absolute value $\left|R_{N} \right| \le R_{N}^{(1)} + R_{N}^{(2)}$:
\begin{align*}
\left|R_{N} \right| \leq N^{\alpha + \alpha \gamma_{f}} W_{B}^{\gamma_f} \int_{0}^{1} \left| c(s) - s \right|^{\delta \gamma_{f}} ds + {N^{\alpha } C_{f } \left| 1 - \tau_{N} \right| }.
\end{align*}
Thus, $\left| R_{N} \right|^{2}  \leq 2 \left( R_{N}^{(1)} \right)^{2} + 2 \left( R_{N}^{(2)} \right)^{2}$. By the Cauchy-Schwarz inequality, the independence between $B^{H_{1}}$ and $\tau_{i}$, and Lemma \ref{lem:step}, recognizing that $\mathbb{E}[W_B^{2\gamma_f}] < \infty$, the expected values are bounded by:
\begin{align}
\mathbb{E} \left[ \left( R_{N}^{(1)} \right)^{2} \right] &\leq \dfrac{C }{N^{\left( \frac{b(\delta \gamma_{f}+1)}{2} - 1 - 2 \alpha \gamma_{f} - 2 \alpha \right) \wedge \left( a(\delta \gamma_{f}+1) - 2 - 2 \alpha \gamma_{f} - 2 \alpha \right) }}, \label{rn1F} \\
\mathbb{E} \left[ \left( R_{N}^{(2)} \right)^{2} \right] &\leq {N^{2 \alpha }} C_{f}^{2} \mathbb{E} \left[ \left| 1 - \tau_{N} \right|^{2} \right] \leq \dfrac{C}{N^{a-1 - 2\alpha}}. \label{rn2F}
\end{align}

\begin{table}[H]
\centering
\begin{tabular}{|l|c|c|}
\hline
 & Relation between $a$ and $\alpha$ & Relation between $b$ and $\alpha$ \\ \hline
$\mathbb{E} \left[ \left( R_{N}^{(1)} \right)^{2} \right]$ & $\alpha < \dfrac{a(H_{1} \gamma_{f} + 1) - 2}{2(\gamma_{f} + 1)}$ & $\alpha < \dfrac{b(H_{1} \gamma_{f} +1) - 2}{4(\gamma_{f} + 1)}$ \\ \hline
$\mathbb{E} \left[ \left( R_{N}^{(2)} \right)^{2} \right]$ & $ \alpha < {\dfrac{a-1}{2}}$ & Not necessary \\ \hline
\end{tabular}
\caption{Order of exponents for $L^{2}$ convergence of $R_N$}
\end{table}

We conclude that $L^{2}$ convergence of $R_{N}$ is ensured if $0 < \alpha < \min \left( \frac{a(H_{1} \gamma_{f} +1)-2}{2(\gamma_{f} +1)}, \frac{a-1}{2} \right)$.

By applying analogous algebraic steps, we decompose the numerator approximation $L_N = L_{N}^{(1)} - L_{N}^{(2)}$, where $L_N^{(2)} = r(x)R_N$, which converges to zero. The term $L_N^{(1)}$ splits into three components: $L_N^{(1.1)}$ and $L_N^{(1.2)}$, dependent on the H\"older continuity of $r(\cdot)$ and $K(\cdot)$, respectively, and $L_N^{(1.3)}$, dependent on the random temporal gap $\left| \tau_N - 1 \right|$. 

For the first term, we explicitly bound its absolute value using the H\"older continuity of $r$ and the boundedness of $K$ (Assumption \ref{assump:regularity}):
\begin{align*}
\left| L_{N}^{(1.1)} \right| &= N^{\alpha} \int_{0}^{1} \left| r \left( B_{c(s)}^{H_{1}} \right) - r \left( B_{s}^{H_{1}} \right) \right| K \left( \left( x - B_{c(s)}^{H_{1}} \right) N^{\alpha}\right) ds \\
&\leq C_K N^{\alpha} \int_{0}^{1} \left| B_{c(s)}^{H_{1}} - B_{s}^{H_{1}} \right|^{\gamma_{r}} ds \le C_K N^{\alpha} W_B^{\gamma_r} \int_{0}^{1} \left| c(s) - s \right|^{\delta \gamma_{r}} ds.
\end{align*}
Squaring and taking the expectation, the independence of $B^{H_1}$ and $\tau_i$ allows us to directly apply Lemma \ref{lem:step}. Since $W_B$ has finite moments of all orders, $\mathbb{E} \left[ W_{B}^{2\gamma_r} \right] < \infty$, yielding:
\begin{align*}
\mathbb{E} \left[ \left| L_{N}^{(1.1)} \right|^{2} \right] &\leq C N^{2 \alpha} \mathbb{E} \left[ W_{B}^{2\gamma_r} \right] \mathbb{E} \left[ \left( \int_{0}^{1} \left| c(s) - s \right|^{\delta \gamma_{r}} ds \right)^{2} \right] \\
&\leq \dfrac{C}{N^{\left[ \frac{b(\delta \gamma_{r}+1)}{2} - 1  - 2 \alpha \right] \wedge \left[ a(\delta \gamma_{r}+1) - 2 - 2 \alpha \right] }}.
\end{align*}

By an analogous algebraic procedure, applying the respective H\"older property for $L_N^{(1.2)}$ and the renewal moment bounds from \eqref{fourth} for $L_N^{(1.3)}$, we obtain the remaining bounds:
\begin{align*}
\mathbb{E} \left[ \left| L_{N}^{(1.2)} \right|^{2} \right] &\leq \dfrac{C_r C_K C_B }{N^{\left[ \frac{b(\delta \gamma_{K}+1)}{2} - 1 - 2\alpha -2 \alpha \gamma_{K} \right] \wedge \left[ a(\delta \gamma_{K}+1) - 2 - 2\alpha -2 \alpha \gamma_{K} \right] }}, \\
\mathbb{E} \left[ \left| L_{N}^{(1.3)} \right|^{2} \right] &\leq \dfrac{C_K C_r}{N^{\left( \frac{2a - 2}{2} - 2 \alpha \right)  \wedge \left( \frac{b - 1}{2} - 2 \alpha \right)}}.
\end{align*}

Before summarizing the final convergence conditions, a formal remark on the H\"older exponent $\delta$ is required.

\begin{rem} \label{rem:delta_limit}
By Assumption \ref{assump:regularity}, the sample paths of the fBm are locally H\"older continuous for any exponent $\delta \in (0, H_1)$, but strictly not for $H_1$. Consequently, the algebraic bounds derived above fundamentally depend on $\delta$. However, notice that the convergence constraints dictate strictly upper bounds on $\alpha$ (e.g., $\alpha < \frac{a(\delta\gamma_r+1)-2}{2(\gamma_r+1)}$). Since these bounding functions are continuous and strictly increasing with respect to $\delta$, if $\alpha$ satisfies the strict inequality evaluated at the supremum limit as $\delta \to H_1^-$, there inherently exists a valid H\"older exponent $\delta < H_1$ sufficiently close to $H_1$ such that the bound remains satisfied. Therefore, replacing $\delta$ with $H_1$ to define the final supremum convergence region is correct.
\end{rem}

Applying the logic of Remark \ref{rem:delta_limit} to the bounds obtained for $L_N^{(1.1)}$, $L_N^{(1.2)}$, and $L_N^{(1.3)}$, we derive the asymptotic constraints summarized in Table \ref{tab:ln_convergence}. \\

\begin{table}[H] 
\centering
\renewcommand{\arraystretch}{1.5}
\begin{tabular}{|l|c|c|}
\hline
 & Relation between $a$ and $\alpha$ & Relation between $b$ and $\alpha$ \\ \hline
$\mathbb{E} \left[ \left| L_{N}^{(1.1)} \right|^{2} \right]$ & $\alpha < \dfrac{a(H_{1} \gamma_{r} +1) - 2}{2} $ & $\alpha < \dfrac{b(H_{1} \gamma_{r} +1) - 2}{4}$ \\ \hline
$\mathbb{E} \left[ \left| L_{N}^{(1.2)} \right|^{2} \right]$ & $\alpha < \dfrac{a(H_{1} \gamma_{K} +1) - 2}{2(1 + \gamma_{K})}$ & $\alpha < \dfrac{b(H_{1} \gamma_{K} +1) - 2}{4(1 + \gamma_{K})}$ \\ \hline
$\mathbb{E} \left[ \left| L_{N}^{(1.3)} \right|^{2} \right]$ & $ \alpha < \dfrac{a-1}{2}$ & $\alpha < \dfrac{b-1}{4}$ \\
\hline
\end{tabular}
\caption{Supremum conditions on $\alpha$ (as $\delta \to H_1^-$) to ensure $L^{2}$-convergence of $L_N$.}
\label{tab:ln_convergence}
\end{table}

We deduce that $L^{2}$ convergence of $L_{N}$ holds if $0 < \alpha < \min \left( \frac{a(H_{1} \gamma_{r} +1)-2}{2}, \frac{a(H_{1} \gamma_{K} +1)-2}{2(\gamma_{K} +1)}, \frac{a-1}{2} \right)$.
\end{proof}
\subsection*{Proof of Lemma \ref{lem:localtime}}
\begin{proof}
Applying the triangle inequality to separate the continuous and discrete integrals over the random intervals gives:
\begin{align}
&\left| N^{\alpha} \int_{0}^{1} f \left( \left( x - B_{c(s)}^{H_1} \right) N^{\alpha} \right) ds - d_{1} L^{H_{1}} (1,x) \right| \nonumber \\
&\leq \left| N^{\alpha} \int_{0}^{1} f \left( \left( x - B_{c(s)}^{H_1} \right) N^{\alpha} \right) ds - N^{\alpha} \int_{0}^{1} f\left( \left( x - B_{s}^{H_{1}} \right) N^{\alpha} \right) ds \right| \label{lt11} \\
&\quad + \left| N^{\alpha} \int_{0}^{1} f\left( \left( x - B_{s}^{H_{1}} \right) N^{\alpha} \right) ds - d_{1} L^{H_{1}} (1,x) \right|. \label{lt12}
\end{align}
The first difference \eqref{lt11} corresponds exactly to the absolute value of the approximation error controlled by $R_N$ in Lemma \ref{lem:integrals}, guaranteeing its $L^2$-convergence to zero. The second difference \eqref{lt12} converges to zero via the deterministic occupation time formula for fBm derived by \cite{sued2013}.

For the numerator term involving the regression function $r(\cdot)$, the procedure is logically analogous. Recalling that the discrete Riemann-Stieltjes sum can be expressed as an integral over the step function $c(s)$, we replace $f(\cdot)$ with the product $r(B_{c(s)}^{H_1})f(\cdot)$ and decompose the absolute error:
\begin{align}
    &\left| N^{\alpha} \int_{0}^{1} r\big(B_{c(s)}^{H_{1}}\big) f\Big(\big(x - B_{c(s)}^{H_{1}}\big)N^{\alpha}\Big) ds - d_{1} r(x) L^{H_{1}}(1,x) \right| \nonumber \\
    &\leq \left| N^{\alpha} \int_{0}^{1} r\big(B_{c(s)}^{H_{1}}\big) f\Big(\big(x - B_{c(s)}^{H_{1}}\big)N^{\alpha}\Big) ds - N^{\alpha} \int_{0}^{1} r\big(B_{s}^{H_{1}}\big) f\Big(\big(x - B_{s}^{H_{1}}\big)N^{\alpha}\Big) ds \right| \label{lt21} \\
    &\quad + \left| N^{\alpha} \int_{0}^{1} r\big(B_{s}^{H_{1}}\big) f\Big(\big(x - B_{s}^{H_{1}}\big)N^{\alpha}\Big) ds - d_{1} r(x) L^{H_{1}}(1,x) \right|. \label{lt22}
\end{align}
Notice that the term in \eqref{lt21} is exactly the absolute value of the approximation error $|L_N|$ defined in Lemma \ref{lem:integrals} (with $f$ playing the role of the kernel). By decomposing $L_N$ into its bounded components, Lemma \ref{lem:integrals} guarantees that \eqref{lt21} converges to $0$ in $L^2$. 

Meanwhile, the term in \eqref{lt22} represents the convergence of the continuous fractional integral. From the results provided in \cite{sued2013} (specifically, Lemma 3.5), and recalling that $d_1 = \int_{\mathbb{R}} f(y) dy$, we can ensure the convergence of \eqref{lt22} to zero. Therefore, the result holds for both the numerator and denominator formulations.
\end{proof}

\subsection*{Proof of Lemma \ref{lem:noise}}
\begin{proof}
The starting point is to examine the second moment of the fractional noise term, taking into account the variance normalization derived in \eqref{constante-m}. By applying the law of total expectation, we expand $\mathbb{E} [ Q_{N}^{2} ]$ into a diagonal component $Q_{N}^{(1.1)}(x)$ and an off-diagonal component $Q_{N}^{(1.2)}(x)$:
\begin{align*}
\mathbb{E} \left[ Q_{N}^{2} \right] &\leq \kappa_{1}^{-2H_{2}} N^{2\alpha + a H_{2}} \sum_{i=0}^{N-1} \mathbb{E} \left[ \Delta \tau_{i}^{2} K^{2} \left( \left( x - B_{\tau_{i}}^{H_{1}} \right) N^{\alpha} \right) \left( B_{\tau_{i+1}}^{H_{2}} - B_{\tau_{i}}^{H_{2}} \right)^{2} \right] \\
&\quad + 2 \kappa_{1}^{-2H_{2}} N^{2\alpha + a H_{2}} \sum_{i < j} \mathbb{E} \left[ \Delta \tau_{i} \Delta \tau_{j} K_{i} K_{j} \left( B_{\tau_{i+1}}^{H_{2}} - B_{\tau_{i}}^{H_{2}} \right) \left( B_{\tau_{j+1}}^{H_{2}} - B_{\tau_{j}}^{H_{2}} \right) \right] \\
&:= Q_{N}^{(1.1)} (x) + Q_{N}^{(1.2)} (x),
\end{align*}
where $K_i$ serves as shorthand for the kernel evaluation at step $i$. We condition the expectations on the independent renewal increments. For $Q_{N}^{(1.1)}$, we condition on $\tau_i = z_i$ and $t_{i+1} = w_1$. For $Q_{N}^{(1.2)}$, we introduce the random gap $X_{j-i-1} = \tau_{j} - \tau_{i+1}$ defined by Araya et al. (2023), conditioning on $X_{j-i-1} = w_2$ and $t_{j+1} = w_3$. 

To obtain a mathematically tractable expression for the off-diagonal covariance, we apply a double Taylor-type linearization to the fractional increments:
\begin{equation} \label{taylor-T22}
\mathbb{E} \left[ \left( B_{z_{i} + w_{1}}^{H_{2}} - B_{z_{i}}^{H_{2}} \right) \left( B_{z_{i} + w_{1} + w_{2} + w_{3}}^{H_{2}} - B_{z_{i} + w_{1} + w_{2}}^{H_{2}} \right) \right] \leq 2H_{2} \left( 2H_{2}-1 \right) w_{1} w_{2}^{2H_{2}-2} w_{3}.
\end{equation}

Crucially, this bound introduces the factor $w_2^{2H_2-2}$. Since $H_2 < 1$, the exponent is strictly negative ($2H_2 - 2 < 0$), creating a potential singularity as $w_2 \to 0$. However, integrating this factor against the density $G_{X_{j-i-1}}(w_2)$ is mathematically licit. The random gap $X_{j-i-1}$ is a convolution of strictly positive inter-arrival times, ensuring its probability density is sufficiently regular near the origin. By explicitly invoking Assumption \ref{assump:inverse} (and the properties discussed in Remark 1), the inverse fractional moments of $X_{j-i-1}$ are guaranteed to be finite, preventing any singularity in the subsequent integration steps.

\noindent \textbf{I. Convergence for Bounded Support Kernels:} \\
Assume the kernel has bounded support on $[-1,1]$. For the diagonal term, applying the Cauchy-Schwarz inequality isolates the indicator functions. Evaluating the bounded integral yields a scale factor of $N^{-\alpha/2}$. Integrating over the renewal densities generates fractional moments:
\begin{align*}
Q_{N}^{(1.1)} (x) &\leq C_k^2 \kappa_{1}^{2H_{2}} N^{2\alpha + a H_{2}} \sum_{i=0}^{N-1} \int_{0}^{\infty} \dfrac{1}{\sqrt[4]{2 \pi z_{i}^{2H_{1}}}} \dfrac{\sqrt{2}}{N^{\alpha/2}} G_{\tau_{i}}(z_{i}) dz_{i} \left( \mathbb{E} \left[ t_{i+1}^{4} \right] \right)^{\frac{H_{2}+1}{2}} \\
&\leq C N^{3\alpha/2 + a H_{2}} \sum_{i=0}^{N-1} \mathbb{E} \left[ \tau_{i}^{-H_{1}/2} \right] \dfrac{\kappa_{2}^{H_{2}/2}}{N^{b (H_{2}+1)/2}} \leq C N^{3\alpha /2 + 1 + a H_{2} - b (H_{2}+1)/2}.
\end{align*}
Since $\tau_i > 0$ almost surely, $\mathbb{E} [ \tau_{i}^{-H_{1}/2} ]$ is strictly bounded under Assumption \ref{assump:inverse}.

For the off-diagonal term $Q_N^{(1.2)}$, substituting \eqref{taylor-T22}, applying Cauchy-Schwarz to separate the kernels, and leveraging the Arithmetic Mean-Geometric Mean (AM-GM) inequality allows us to completely decouple the variables in the denominator:
\begin{align*}
Q_{N}^{(1.2)} (x) &\leq C N^{3 \alpha /2 - 2 +a H_{2}} \sum_{i < j} \int_{0}^{\infty} \dots \int_{0}^{\infty} \dfrac{w_{1} w_{2}^{2H_{2}-2} w_{3}}{\sqrt[8]{ z_{i}^{8H_{1}/3} w_{1}^{2H_{1}/3} w_{2}^{2H_{1}/3} }} G_{\tau_{i}} g_{t_{i+1}} G_{X_{j-i-1}} g_{t_{j+1}} d\mathbf{w} \\
&\leq C N^{3 \alpha /2 - 2 +a H_{2}} \sum_{j=1}^{N-1} \sum_{i=0}^{j-1} \mathbb{E} \left[ \tau_{i}^{-H_{1}/3} \right] \mathbb{E} \left[ t_{i+1}^{1 - H_{1}/12} \right] \mathbb{E} \left[ X_{j-i-1}^{-(2-2H_{2}+H_{1}/12)} \right] \mathbb{E} \left[ t_{j+1} \right].
\end{align*}
The positivity of the exponents verifies Assumption \ref{assump:inverse} for the random times $\tau_i$ and $X_{j-i-1}$. Evaluating the sums yields the upper bound $O(N^{3 \alpha /2 - 2 + H_{1}/12 + aH_{2}})$. Combining both diagonal and off-diagonal bounds, $L^2$ convergence to $0$ is guaranteed if:
$$ \alpha < \min \left( \frac{2} {3}\left( \frac{b}{2} - 1 - H_{2} \left( a - \frac{b}{2} \right) \right) , \frac{2}{3} \left( 2 - \frac{H_{1}}{12} - aH_{2} \right)\right). $$

\noindent \textbf{II. Convergence for the Gaussian Kernel:} \\
For the Gaussian kernel, we explicitly integrate over the joint normal distribution. For the diagonal term, standard Gaussian integration yields:
\begin{align*}
Q_{N}^{(1.1)} (x) &= \dfrac{\kappa_{1}^{-2H_{2}}}{2 \pi^{3/2}} N^{2\alpha + a H_{2}} \sum_{i=0}^{N-1} \int_{0}^{\infty} \dfrac{1}{\sqrt{2N^{2 \alpha} z_{i}^{2H_{1}} +1}} G_{\tau_{i}}(z_{i}) dz_{i} \int_{0}^{\infty} w_{1}^{2H_{2}+2} g_{t_{i+1}}(w_{1}) dw_{1} \\
&\leq C_{H_1} N^{\alpha + a H_{2}} \sum_{i=0}^{N-1} \mathbb{E} \left[ \tau_i^{-H_1} \right] \mathbb{E} \left[ \left( t_{i+1}^{4} \right)^{\frac{H_2 + 1}{2}} \right] \leq C N^{\alpha + a H_{2} + 1 - b\left( \frac{H_2 +1}{2} \right)},
\end{align*}
where Jensen's inequality was applied to the fractional moment of $t_{i+1}$ (Assumption \ref{assump:renewal}).

For the off-diagonal term, the expectation requires integrating the bi-variate Gaussian density of $\widetilde{X} = ( B_{z_{i}}^{H_{1}}, B_{z_{i} + w_{1} + w_{2}}^{H_{1}} )$. By applying the deterministic lower bound for the fractional covariance determinant, $\det(\Sigma) \geq \kappa \left( w_{1} + w_{2} \right)^{2H_{1}} z_{i}^{2H_{1}}$ (\cite{yan2009}), the joint integral is bounded:
\begin{align*}
Q_{N}^{(1.2)} (x) &\leq C N^{2 \alpha +a H_{2}} \sum_{i < j} \int \dots \int \dfrac{ w_{1}^{2} w_{2}^{2H_{2}-2} w_{3}^{2}}{\sqrt{\det(\Sigma) N^{4 \alpha} + N^{2 \alpha} z_{i}^{2H_{1}} + N^{2 \alpha} (z_{i} + w_{1} + w_{2})^{2H_{1}} + 1}} d\mathbf{w} \\
&\leq C N^{\alpha - 2 +aH_{2}} \sum_{i < j} \mathbb{E} \left[ \tau_{i}^{-7H_{1}/12} \right] \mathbb{E} \left[ t_{i+1}^{2 - 5H_{1}/24} \right] \mathbb{E} \left[ X_{j-i-1}^{-(2 - 2H_{2} + 5H_{1}/24)} \right] \mathbb{E} \left[ t_{j+1}^{2} \right].
\end{align*}
Using Assumption \ref{assump:inverse} for the random gaps, and applying Jensen's inequality to $t_{i+1}$, the nested double sum collapses to $N^2$, giving:
\begin{align*}
Q_{N}^{(1.2)} (x) &\leq C N^{\alpha +aH_{2} -a(1 -5H_{1}/48) - a} \cdot N^{2} = C N^{\alpha -a \left( 2 - H_2 - \frac{5 H_1}{48} \right) +2}.
\end{align*}
Combining both constraints, the variance of the Gaussian fractional noise strictly converges to $0$ if $\alpha < \min \left( \frac{b}{2} - 1 - H_2 \left( a - \frac{b}{2} \right), a\left( 2 - H_{2} - \frac{5H_{1}}{48} \right) - 2 \right)$.
\end{proof}
\subsection*{Proof of Proposition \ref{prop:fixedN}}
\begin{proof}
The $L^2$-consistency for the fixed sample size $N$ is established by partitioning the estimation error into a noise-free approximation and a fractional perturbation:
\begin{equation*}
\hat{r}_N(x) - r(x) = \left[ T_1^{(N)}(x) - r(x) \right] + T_2^{(N)}(x).
\end{equation*}

From Lemma \ref{lem:integrals}, the discrete Riemann-Stieltjes sums in $T_1^{(N)}(x)$ can be replaced by continuous integrals over $[0, \tau_N]$ with an asymptotically vanishing error. By the Law of Large Numbers (Assumption \ref{assump:renewal}), the final random time $\tau_N \to 1$ almost surely. This allows the direct application of the occupation time formula derived by \cite{sued2013} (Lemma \ref{lem:localtime}), ensuring that the numerator and denominator of $T_1^{(N)}(x)$ converge to proportional fractional limits that cancel out exactly to $r(x)$. Thus, $T_1^{(N)}(x) \xrightarrow{L^2} r(x)$. 

Simultaneously, Lemma \ref{lem:noise} explicitly controls the variance of the long-memory noise perturbation, guaranteeing $T_2^{(N)}(x) \xrightarrow{L^2} 0$. The intersection of the necessary parametric constraints for each component is summarized in Table \ref{tab:convergence_fixed_N}.
\end{proof}

\begin{table}[H]
\centering
\renewcommand{\arraystretch}{1.5}
\begin{tabular}{|l|l|}
\hline
\textbf{Component analyzed} & \textbf{$L^{2}$-convergence constraint on $\alpha$} \\ \hline
$\mathbb{E} \left[ \left( R_{N} \right)^{2} \right]$ (Lemma \ref{lem:integrals}) & 
    $\alpha < \min \left( \frac{a(H_{1} \gamma_{f} + 1) - 2}{2(\gamma_{f} + 1)}, \, \frac{a-1}{2} \right)$ \\ \hline
$\mathbb{E} \left[ \left( L_{N} \right)^{2} \right]$ (Lemma \ref{lem:integrals}) & 
    $\alpha < \min \left( \frac{a(H_{1} \gamma_{r} + 1) - 2}{2(\gamma_{r} + 1)}, \, \frac{a(H_{1} \gamma_{K} + 1) - 2}{2(\gamma_{K} + 1)}, \, \frac{a-1}{2(1 + \gamma_{K})} \right)$ \\ \hline
$\mathbb{E} \left[ Q_{N}^{2} \right]$ - Bounded Kernel (Lemma \ref{lem:noise}) & 
    $\alpha < \min \left( \frac{2}{3}\left( \frac{b}{2} - 1 - H_{2} \left( a - \frac{b}{2} \right) \right), \, \frac{2}{3} \left( 2 - \frac{H_{1}}{12} - aH_{2} \right) \right)$ \\ \hline
$\mathbb{E} \left[ Q_{N}^{2} \right]$ - Gaussian Kernel (Lemma \ref{lem:noise}) & 
    $\alpha < a\left( 2 - H_{2} - \frac{5H_{1}}{48} \right) - 2$ \\ \hline
\end{tabular}
\caption{Upper bounds for the bandwidth parameter $\alpha$ to ensure $L^2$-convergence under a fixed sample size $N$ (Proposition \ref{prop:fixedN}).}
\label{tab:convergence_fixed_N}
\end{table}

\subsection*{Proof of Proposition \ref{prop:randomN}}
\begin{proof}
To control the discrepancy introduced by the random sampling size $N(1)$, we partition the probability space. Define the subsets governing the counting process, the stopping time, and the denominators:
\begin{align*}
\Omega_N^1 &= \left\{ \left| \frac{N(1)}{N} - 1 \right| \leq N^{-\delta} \right\}, \quad \mathbb{P}(\Omega_N^1) \geq 1 - \frac{C_1}{N^\mu} \quad \text{with } \mu = \min(2a-2-2\delta,\, b-1-3\delta), \\
\Omega_N^2 &= \left\{ |1 - \tau_N| \leq N^{-\gamma} \right\}, \quad \mathbb{P}(\Omega_N^2) \geq 1 - \frac{C_2}{N^\nu} \quad \text{with } \nu = -4\gamma + \min(2a-2,\, b-1), \\
\Omega_N^3 &= \{ D_N(x) \geq c \;\wedge\; D_{N(1)}(x) \geq c \}.
\end{align*}
Given the $L^2$-convergence to a strictly positive local time ($L^{H_1}(1,x) > 0$ a.s.), $\mathbb{P}(\Omega_N^3) \geq 1 - C_3/N^\eta$. Defining the regular set $\Omega_N = \Omega_N^1 \cap \Omega_N^2 \cap \Omega_N^3$, Boole’s inequality yields $\mathbb{P}(\Omega_N^c) \leq C_1 N^{-\mu} + C_2 N^{-\nu} + C_3 N^{-\eta}$.

We partition the expected squared difference of the estimators:
\begin{equation} \label{eq:partition_error}
\mathbb{E} \left[ \left| \widehat{r}_{N} - \widehat{r}_{N(1)} \right|^{2} \right] = \mathbb{E} \left[ \left| \widehat{r}_{N} - \widehat{r}_{N(1)} \right|^{2} \mathbbm{1}_{\Omega_{N}} \right] + \mathbb{E} \left[ \left| \widehat{r}_{N} - \widehat{r}_{N(1)} \right|^{2} \mathbbm{1}_{\Omega_{N}^{c}} \right].
\end{equation}

On the complement set $\Omega_N^c$, applying the Cauchy-Schwarz and $C_p$ inequalities yields:
\begin{align*}
\mathbb{E} \left[ \left| \widehat{r}_{N} - \widehat{r}_{N(1)} \right|^{2} \mathbbm{1}_{\Omega_{N}^{c}} \right] &\leq \sqrt{\mathbb{E} \left[ \left| \widehat{r}_{N} - \widehat{r}_{N(1)} \right|^{4} \right]} \sqrt{\mathbb{P} \left( \Omega_{N}^{c} \right)} \\
&\leq \sqrt{8 \left( \mathbb{E}\left[ \left| \widehat{r}_{N} \right|^4\right] + \mathbb{E}\left[\left| \widehat{r}_{N(1)} \right|^4 \right] \right)} \sqrt{\mathbb{P} \left( \Omega_{N}^{c} \right)}. 
\end{align*}
Because the numerators are driven by the noise process $\varepsilon^{H_2}$ over a compact interval, Fernique's Theorem guarantees that the fractional Gaussian paths possess bounded marginal moments ($\sup_{t} \mathbb{E}[|\varepsilon_t^{H_2}|^4] \le M_4 < \infty$). Consequently, the expectation on $\Omega_N^c$ vanishes if the failure probabilities decay ($\mu > 0, \nu > 0$), which structurally imposes $\alpha < \gamma < \frac{1}{4} \min (2a - 2, b - 1)$.

On the regular set $\Omega_N$, denominators are deterministically bounded below by $c$. Expanding $\hat{r}_{N(1)} - \hat{r}_N$ via their respective numerators ($A$) and denominators ($D$), we obtain:

\begin{align*}
\left| \widehat{r}_{N} - \widehat{r}_{N(1)} \right|^{2} \mathbbm{1}_{\Omega_{N}} &\leq \left[ \dfrac{2}{c^2} \left( A_{N(1)} - A_{N} \right)^{2} + \dfrac{2}{c^4} A_{N}^{2} \left( D_{N} - D_{N(1)} \right)^{2} \right] \mathbbm{1}_{\Omega_{N}} \\
&\leq C N^{2\alpha} \left| \tau_{N(1)} - \tau_N \right|^2 \mathbbm{1}_{\Omega_{N}}.
\end{align*}
To control this temporal gap, the triangle inequality establishes $\left| \tau_{N(1)} - \tau_N \right| \leq \left| \tau_{N(1)} - 1 \right| + \left| \tau_{N} - 1 \right|$. By the counting process definition, $\tau_{N(1)} \leq 1 < \tau_{N(1)+1}$, which implies $0 \leq 1 - \tau_{N(1)} < \Delta \tau_{N(1)}$. Since the inter-arrival gaps vanish at $\mathcal{O}(N^{-1})$ and $|\tau_N - 1| \leq N^{-\gamma}$ on $\Omega_N^2$, the total gap is bounded by $
\left| \tau_{N(1)} - \tau_N \right| \mathbbm{1}_{\Omega_{N}} \leq C_T N^{-\gamma} \mathbbm{1}_{\Omega_{N}}$. 

Selecting $\alpha < \gamma$ guarantees that $\mathbb{E} \left[ \left| \widehat{r}_{N} - \widehat{r}_{N(1)} \right|^{2} \mathbbm{1}_{\Omega_{N}} \right] \to 0$. Both expectation terms in \eqref{eq:partition_error} successfully collapse to zero, completing the proof.
\end{proof}
\subsection*{Proof of Theorem \ref{thm:main}}
\begin{proof}
To establish the $L^2$-consistency of the wN-W estimator $\hat{r}_{N(1)}(x)$ under a random number of observations, we decompose the estimation error. By applying the Cauchy-Schwarz inequality, the overall mean squared error is bounded by:
\begin{equation} \label{eq:theorem_decomposition}
\mathbb{E} \left[ \big| \hat{r}_{N(1)}(x) - r(x) \big|^2 \right] \le 2  \mathbb{E} \left[ \big| \hat{r}_{N(1)}(x) - \hat{r}_N(x) \big|^2 \right] + 2  \mathbb{E} \left[ \big| \hat{r}_N(x) - r(x) \big|^2 \right].
\end{equation}

The right-most term of \eqref{eq:theorem_decomposition} represents the fixed-sample estimation error, which converges to zero as $N \to \infty$ provided the bandwidth parameter $\alpha$ satisfies the constraints established in Proposition \ref{prop:fixedN}. 

The first term on the right-hand side represents the stochastic discrepancy induced by the random sampling times $\tau_i$. As demonstrated in Proposition \ref{prop:randomN}, the expectation over the complement set $\Omega_N^c$ decays to zero given the moment constraints $a>1$, $b>1$, and $\alpha < \frac{1}{4} \min(2a-2, b-1)$. Simultaneously, on the regular set $\Omega_N$, the strictly bounded temporal gap $|\tau_{N(1)} - \tau_N|$ guarantees that the conditional expectation also collapses to zero.

Therefore, global $L^2$-consistency is strictly guaranteed by satisfying the intersection of all admissible parametric sets derived from both propositions. These intersecting bounds correspond exactly to the conditions formulated in Theorem \ref{thm:main}, completing the proof.
\end{proof} 

To understand the practical implications of the theoretical bounds in Theorem \ref{thm:main}, Table \ref{tabla-alpha} provides an overview of the maximum admissible values for the bandwidth parameter $\alpha$ under various degrees of temporal irregularity ($a, b$) and noise persistence ($H_2$), keeping the covariate memory fixed at $H_1 = 0.55$. 

\begin{table}[H]
\centering
\begin{tabular}{llcccccc}
\hline
\multirow{2}{*}{\textbf{$a$}} & \multirow{2}{*}{\textbf{$b$}} & \multicolumn{3}{c}{\textbf{Quartic Kernel ($H_2$)}} & \multicolumn{3}{c}{\textbf{Gaussian Kernel ($H_2$)}} \\ \cline{3-8}
 & & \textbf{0.60} & \textbf{0.75} & \textbf{0.90} & \textbf{0.60} & \textbf{0.75} & \textbf{0.90} \\ \hline
1.60 & 3.4 & \multicolumn{6}{c}{ Not valid ($b \ge 2a$)} \\ \hline
1.75 & 1.5 & 0.1250 & 0.1250 & 0.0667 & 0.1250 & 0.1250 & 0.1250 \\
1.75 & 2.5, 3.4 & 0.1781 & 0.1781 & 0.1781 & 0.1781 & 0.1781 & 0.1781 \\ \hline
1.85 & 1.5 & 0.1250 & 0.1167 & 0.0067 & 0.1250 & 0.1250 & 0.1250 \\
1.85 & 2.5, 3.4 & 0.2125 & 0.2125 & 0.1928 & 0.2125 & 0.2125 & 0.2125 \\ \hline
1.95 & 1.5 & 0.1250 & 0.0667 & { Invalid} & 0.1250 & 0.1250 & 0.1250 \\
1.95 & 2.5, 3.4 & 0.2375 & 0.2375 & 0.1328 & 0.2375 & 0.2375 & 0.1304 \\ \hline
\end{tabular}%
\caption{Admissible upper bounds for $\alpha$ under selected parameter combinations ($H_1 = 0.55$).}
\label{tabla-alpha}
\end{table}

As detailed in Table \ref{tabla-alpha}, the interaction between highly irregular sampling times and strongly persistent noise can significantly reduce the admissible region for $\alpha$. In extreme scenarios, this combination leads to theoretically invalid configurations, either by violating the inequality for the renewal process ($b \ge 2a$) or by forcing a negative bandwidth exponent ($\alpha \le 0$), making $L^2$ consistency impossible.

To facilitate an interactive exploration of the feasible region of the smoothing parameter in relation to various long-memory settings and renewal times, an open-source web application has been developed and is available at: \url{https://f78s54-tania0roa0rojas.shinyapps.io/ValidSubset/}. This graphical tool allows readers and practitioners to dynamically adjust system parameters and visualize the resulting set of admissible values for $\alpha$, automatically identifying in real-time which mathematical component acts as the dominant constraint. This application serves as a bridge between our theoretical findings and their practical implementation in subsequent simulation studies.

\section{Simulation Results} \label{simulations}

In this section, we present simulation studies to evaluate the performance of the proposed wN-W estimator $\widehat{r}(x)$ and the underlying stochastic process $Y_{\tau_i}$. We analyze the Empirical Mean Squared Error (MSE) across different expected sample sizes $N$ and memory parameters $H_1, H_2$. For baseline comparisons under regular uniform sampling, we refer to the deterministic framework established by \cite{sued2013}.

For the covariate process driven by fBm, we fix the Hurst parameter at $H_1 = 0.55$. This value approaches the standard Brownian motion boundary ($H_1 = 0.5$), representing a highly variable and unstable trajectory scenario. For the independent fractional noise, we explore three distinct levels of long-memory persistence: $H_2 \in \{0.60, 0.75, 0.90\}$.

The true regression curve is defined by the following ``bump function'' (previously analyzed by \cite{fermin2017}), which is locally H\"older continuous with exponent $\gamma_r = 1$:
\begin{equation} \label{r1}
    r_{1}(x) =  0.7x + 2e^{-10x^{2}}.
\end{equation}

\subsection{Interactive Visualization and Estimation Performance}

The theoretical framework developed in this study involves a complex multidimensional interplay among the long-memory dependence of the noise ($H_2$), the memory of the covariate ($H_1$), the irregularity of the random sampling times ($N(1)$), and the bandwidth parameter ($\alpha$). Since these intricate relationships are difficult to fully capture through static plots, we developed a dynamic and interactive representation.

To provide a comprehensive empirical validation of our theoretical results, we built an open-source web application using the Shiny framework for R, publicly available at \url{https://f78s54-tania0roa0rojas.shinyapps.io/wN-W_estimator/}. 
\begin{figure}[H]
    \centering
    \includegraphics[width=\linewidth]{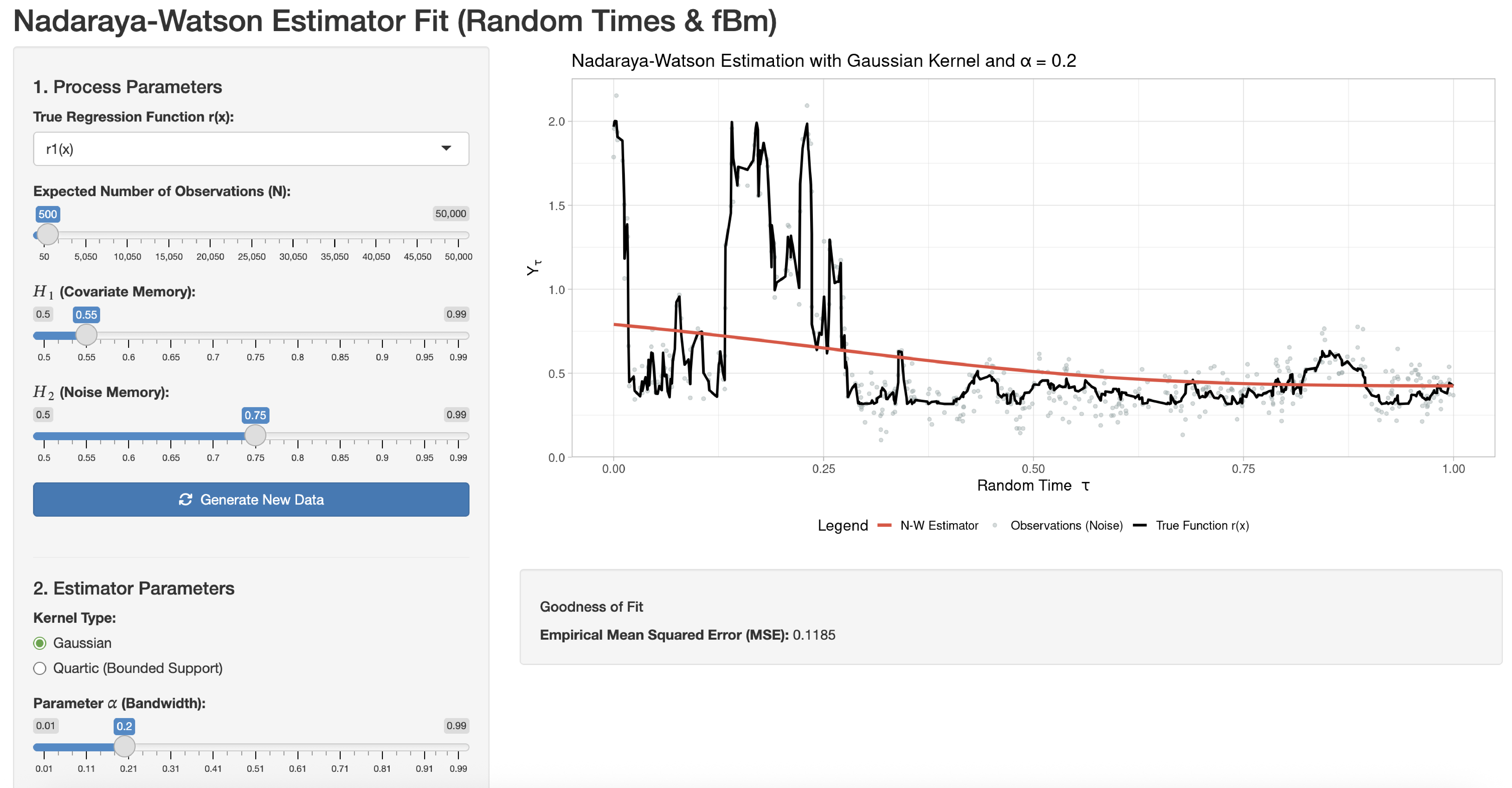}
    \caption{Interactive web application interface for dynamic wN-W estimation.}
    \label{shinyapp}
\end{figure}
This interactive tool empowers readers and practitioners to dynamically generate and analyze the stochastic process $Y_{\tau_i}$ across a continuous spectrum of scenarios. Users can specify the true regression function $r(x)$ and the expected sample size $N$, dynamically simulating the random times $\tau_i$ alongside the fractional Brownian motions for both the covariate ($H_1$) and the noise ($H_2$). Furthermore, the application allows for dynamic estimation by toggling between Gaussian and bounded support kernels while adjusting the bandwidth parameter $\alpha \in (0, 1)$. The interface instantly overlays the resulting wN-W estimator onto the noisy observations and computes the Empirical MSE in real time. By manipulating these parameters, users can visually and quantitatively verify the theoretical bounds established in Theorem 1, directly observing the bias-variance tradeoff and the asymptotic stabilization of the estimator around the true regression function as the irregular temporal gaps shrink for large $N$. \\

To clarify the practical utility of the accompanying Shiny application and address potential concerns regarding methodological circularity, it is crucial to delineate its role as an a priori simulation tool. The interactive application is primarily designed to assist researchers in experimental design and sensitivity analysis, allowing them to explore how different theoretical parameter values ($a$ and $b$) governing the random times impact the asymptotic constraints before data collection. However, in real-world applications—as explicitly demonstrated in the empirical analysis of the Vostok dataset—these parameters are not presumed known. Instead, they are strictly data-driven, estimated directly from the observed temporal gaps ($\Delta \tau_i$) via the method of moments. Consequently, the Shiny application serves as a complementary theoretical sandbox for understanding complex sampling topologies, whereas the practical implementation of the estimator remains entirely empirical and non-circular.

\subsection{Empirical Mean Squared Error Analysis}

To quantitatively evaluate the consistency of the wN-W estimator, we compute the Empirical Mean Squared Error (MSE) over the random observation points. For both kernel functions, the MSE is defined as:
\begin{equation}
MSE \left(\widehat{Y}_{\tau} \right) =  \dfrac{1}{N(1)} \sum_{i=1}^{N(1)} \left( \widehat{Y}_{\tau_{i}}(x) - Y_{\tau_{i}}(x) \right)^{2},  \label{mse}
\end{equation}
where $N(1)$ represents the random number of observations falling within the interval $[0,1]$. By the asymptotic properties of the renewal process established in Assumption \ref{assump:renewal}, the ratio $N(1)/N \to 1$ almost surely, bridging the empirical random count with the deterministic expected sample size $N$. 

Table \ref{mse_tabla} summarizes the computed MSE for the nonparametric estimators under deterministic and random sampling designs, contrasting the two kernel types across different sample sizes ($N$), bandwidth exponents ($\alpha$), and noise persistence levels ($H_2$).
\begin{table}[H]  
\centering
\resizebox{\textwidth}{!}{%
\begin{tabular}{|llllll|llllll|}
\hline
\multicolumn{6}{|c|}{\textbf{Deterministic}}                                                                                                                                                                                                           & \multicolumn{6}{c|}{\textbf{Random}}                                                                                                                                                                                                                 \\ \hline
\multicolumn{1}{|l|}{\textbf{Kernel}}             & \multicolumn{1}{l|}{\textbf{$\alpha$}}     & \multicolumn{1}{l|}{\textbf{N}}                    & \multicolumn{1}{l|}{\textbf{$H_2=0.60$}}               & \multicolumn{1}{l|}{\textbf{$H_2=0.75$}}               & \textbf{$H_2=0.90$}               & \multicolumn{1}{l|}{\textbf{Kernel}}             & \multicolumn{1}{l|}{\textbf{$\alpha$}}     & \multicolumn{1}{l|}{\textbf{N(1)}}                 & \multicolumn{1}{l|}{\textbf{$H_2=0.60$}}                & \multicolumn{1}{l|}{\textbf{$H_2=0.75$}}                & \textbf{$H_2=0.90$}               \\ \hline
\multicolumn{1}{|l|}{}                            & \multicolumn{1}{l|}{}                      & \multicolumn{1}{l|}{50}                            & \multicolumn{1}{l|}{0.2661534}                         & \multicolumn{1}{l|}{0.1516746}                         & 0.1886984                         & \multicolumn{1}{l|}{}                            & \multicolumn{1}{l|}{}                      & \multicolumn{1}{l|}{50}                            & \multicolumn{1}{l|}{0.3590721}                          & \multicolumn{1}{l|}{0.3673873}                          & 0.2937165                         \\ \cline{3-6} \cline{9-12} 
\multicolumn{1}{|l|}{}                            & \multicolumn{1}{l|}{}                      & \multicolumn{1}{l|}{1e+03}                         & \multicolumn{1}{l|}{0.109745}                          & \multicolumn{1}{l|}{0.0947708}                         & 0.0936384                         & \multicolumn{1}{l|}{}                            & \multicolumn{1}{l|}{}                      & \multicolumn{1}{l|}{1e+03}                         & \multicolumn{1}{l|}{0.3213928}                          & \multicolumn{1}{l|}{0.3199888}                          & 0.3183443                         \\ \cline{3-6} \cline{9-12} 
\multicolumn{1}{|l|}{}                            & \multicolumn{1}{l|}{\multirow{-3}{*}{0.10}} & \multicolumn{1}{l|}{\cellcolor[HTML]{FD6864}1e+05} & \multicolumn{1}{l|}{\cellcolor[HTML]{FD6864}0.5738293} & \multicolumn{1}{l|}{\cellcolor[HTML]{FD6864}0.5745912} & \cellcolor[HTML]{FD6864}0.5736525 & \multicolumn{1}{l|}{}                            & \multicolumn{1}{l|}{\multirow{-3}{*}{0.1}} & \multicolumn{1}{l|}{\cellcolor[HTML]{FD6864}1e+05} & \multicolumn{1}{l|}{\cellcolor[HTML]{FD6864}0.4863244}  & \multicolumn{1}{l|}{\cellcolor[HTML]{FD6864}0.4855444}  & \cellcolor[HTML]{FD6864}0.4856166  \\ \cline{2-6} \cline{8-12} 
\multicolumn{1}{|l|}{}                            & \multicolumn{1}{l|}{}                      & \multicolumn{1}{l|}{50}                            & \multicolumn{1}{l|}{0.2619626}                         & \multicolumn{1}{l|}{0.1506687}                         & 0.1768822                         & \multicolumn{1}{l|}{}                            & \multicolumn{1}{l|}{}                      & \multicolumn{1}{l|}{50}                            & \multicolumn{1}{l|}{0.352546}                           & \multicolumn{1}{l|}{0.3553474}                          & 0.2814558                         \\ \cline{3-6} \cline{9-12} 
\multicolumn{1}{|l|}{}                            & \multicolumn{1}{l|}{}                      & \multicolumn{1}{l|}{1e+03}                         & \multicolumn{1}{l|}{0.09192583}                        & \multicolumn{1}{l|}{0.0776438}                         & 0.07483429                        & \multicolumn{1}{l|}{}                            & \multicolumn{1}{l|}{}                      & \multicolumn{1}{l|}{1e+03}                         & \multicolumn{1}{l|}{0.2197067}                          & \multicolumn{1}{l|}{0.2179231}                          & 0.2178756                         \\ \cline{3-6} \cline{9-12} 
\multicolumn{1}{|l|}{}                            & \multicolumn{1}{l|}{\multirow{-3}{*}{0.20}}    & \multicolumn{1}{l|}{1e+05}                         & \multicolumn{1}{l|}{0.3189735}                         & \multicolumn{1}{l|}{0.3189604}                         & 0.3184624                         & \multicolumn{1}{l|}{}                            & \multicolumn{1}{l|}{\multirow{-3}{*}{0.20}} & \multicolumn{1}{l|}{1e+05}                         & \multicolumn{1}{l|}{0.2277019}                          & \multicolumn{1}{l|}{0.2267811}                          & 0.2267169                         \\ \cline{2-6} \cline{8-12} 
\multicolumn{1}{|l|}{}                            & \multicolumn{1}{l|}{}                      & \multicolumn{1}{l|}{50}                            & \multicolumn{1}{l|}{0.2571986}                         & \multicolumn{1}{l|}{0.1494979}                         & 0.1612631                         & \multicolumn{1}{l|}{}                            & \multicolumn{1}{l|}{}                      & \multicolumn{1}{l|}{50}                            & \multicolumn{1}{l|}{0.3425697}                          & \multicolumn{1}{l|}{0.3368326}                          & 0.2644141                         \\ \cline{3-6} \cline{9-12} 
\multicolumn{1}{|l|}{}                            & \multicolumn{1}{l|}{}                      & \multicolumn{1}{l|}{1e+03}                         & \multicolumn{1}{l|}{0.0729148}                         & \multicolumn{1}{l|}{0.0592166}                         & 0.0557591                         & \multicolumn{1}{l|}{}                            & \multicolumn{1}{l|}{}                      & \multicolumn{1}{l|}{1e+03}                         & \multicolumn{1}{l|}{0.1713224}                          & \multicolumn{1}{l|}{0.168998}                           & 0.169112                          \\ \cline{3-6} \cline{9-12} 
\multicolumn{1}{|l|}{\multirow{-9}{*}{Gaussian}} & \multicolumn{1}{l|}{\multirow{-3}{*}{0.25}} & \multicolumn{1}{l|}{1e+05}                         & \multicolumn{1}{l|}{0.1269318}                         & \multicolumn{1}{l|}{0.1264641}                         & 0.1260219                         & \multicolumn{1}{l|}{\multirow{-9}{*}{Gaussian}} & \multicolumn{1}{l|}{\multirow{-3}{*}{0.3}} & \multicolumn{1}{l|}{1e+05}                         & \multicolumn{1}{l|}{0.05788913}                         & \multicolumn{1}{l|}{0.0570189}                          & 0.0567967                         \\ \hline
\multicolumn{1}{|l|}{}                            & \multicolumn{1}{l|}{}                      & \multicolumn{1}{l|}{50}                            & \multicolumn{1}{l|}{0.2565114}                         & \multicolumn{1}{l|}{0.1497628}                         & 0.1578216                         & \multicolumn{1}{l|}{}                            & \multicolumn{1}{l|}{}                      & \multicolumn{1}{l|}{50}                            & \multicolumn{1}{l|}{0.3431026}                          & \multicolumn{1}{l|}{0.3351815}                          & 0.2606817                         \\ \cline{3-6} \cline{9-12} 
\multicolumn{1}{|l|}{}                            & \multicolumn{1}{l|}{}                      & \multicolumn{1}{l|}{1e+03}                         & \multicolumn{1}{l|}{0.0862823}                         & \multicolumn{1}{l|}{0.0721686}                         & 0.0689495                         & \multicolumn{1}{l|}{}                            & \multicolumn{1}{l|}{}                      & \multicolumn{1}{l|}{1e+03}                         & \multicolumn{1}{l|}{0.2003309}                          & \multicolumn{1}{l|}{0.1984449}                          & 0.1985662                         \\ \cline{3-6} \cline{9-12} 
\multicolumn{1}{|l|}{}                            & \multicolumn{1}{l|}{\multirow{-3}{*}{0.10}} & \multicolumn{1}{l|}{1e+05}                         & \multicolumn{1}{l|}{0.3571319}                         & \multicolumn{1}{l|}{0.3571672}                         & 0.3566951                         & \multicolumn{1}{l|}{}                            & \multicolumn{1}{l|}{\multirow{-3}{*}{0.10}} & \multicolumn{1}{l|}{1e+05}                         & \multicolumn{1}{l|}{0.2950395}                          & \multicolumn{1}{l|}{0.2940885}                          & 0.2941242                         \\ \cline{2-6} \cline{8-12} 
\multicolumn{1}{|l|}{}                            & \multicolumn{1}{l|}{}                      & \multicolumn{1}{l|}{50}                            & \multicolumn{1}{l|}{0.2540677}                         & \multicolumn{1}{l|}{0.1485755}                         & 0.1458834                         & \multicolumn{1}{l|}{}                            & \multicolumn{1}{l|}{}                      & \multicolumn{1}{l|}{50}                            & \multicolumn{1}{l|}{0.3256095}                          & \multicolumn{1}{l|}{0.306079}                           & 0.2436255                         \\ \cline{3-6} \cline{9-12} 
\multicolumn{1}{|l|}{}                            & \multicolumn{1}{l|}{}                      & \multicolumn{1}{l|}{1e+03}                         & \multicolumn{1}{l|}{0.0688117}                         & \multicolumn{1}{l|}{0.0553279}                         & 0.0517929                         & \multicolumn{1}{l|}{}                            & \multicolumn{1}{l|}{}                      & \multicolumn{1}{l|}{1e+03}                         & \multicolumn{1}{l|}{0.1643763}                          & \multicolumn{1}{l|}{0.161773}                           & 0.1619564                         \\ \cline{3-6} \cline{9-12} 
\multicolumn{1}{|l|}{}                            & \multicolumn{1}{l|}{\multirow{-3}{*}{0.20}} & \multicolumn{1}{l|}{1e+05}                         & \multicolumn{1}{l|}{0.1608186}                         & \multicolumn{1}{l|}{0.1604468}                         & 0.1599897                         & \multicolumn{1}{l|}{}                            & \multicolumn{1}{l|}{\multirow{-3}{*}{0.20}} & \multicolumn{1}{l|}{1e+05}                         & \multicolumn{1}{l|}{0.0791516}                          & \multicolumn{1}{l|}{0.0782922}                          & 0.0780666                         \\ \cline{2-6} \cline{8-12} 
\multicolumn{1}{|l|}{}                            & \multicolumn{1}{l|}{}                      & \multicolumn{1}{l|}{50}                            & \multicolumn{1}{l|}{0.2421178}                         & \multicolumn{1}{l|}{0.1416007}                         & 0.1353937                         & \multicolumn{1}{l|}{}                            & \multicolumn{1}{l|}{}                      & \multicolumn{1}{l|}{50}                            & \multicolumn{1}{l|}{0.268778}                           & \multicolumn{1}{l|}{0.2409065}                          & 0.1971679                         \\ \cline{3-6} \cline{9-12} 
\multicolumn{1}{|l|}{}                            & \multicolumn{1}{l|}{}                      & \multicolumn{1}{l|}{1e+03}                         & \multicolumn{1}{l|}{0.0474537}                         & \multicolumn{1}{l|}{0.0350902}                         & 0.0324952                         & \multicolumn{1}{l|}{}                            & \multicolumn{1}{l|}{}                      & \multicolumn{1}{l|}{1e+03}                         & \multicolumn{1}{l|}{0.0861950}                          & \multicolumn{1}{l|}{0.0827619}                          & 0.0822441                         \\ \cline{3-6} \cline{9-12} 
\multicolumn{1}{|l|}{\multirow{-9}{*}{Quartic}}  & \multicolumn{1}{l|}{\multirow{-3}{*}{0.25}} & \multicolumn{1}{l|}{\cellcolor[HTML]{67FD9A}1e+05} & \multicolumn{1}{l|}{\cellcolor[HTML]{67FD9A}0.0646294} & \multicolumn{1}{l|}{\cellcolor[HTML]{67FD9A}0.0639921} & \cellcolor[HTML]{67FD9A}0.0636587 & \multicolumn{1}{l|}{\multirow{-9}{*}{Quartic}} & \multicolumn{1}{l|}{\multirow{-3}{*}{0.25}} & \multicolumn{1}{l|}{\cellcolor[HTML]{67FD9A}1e+05} & \multicolumn{1}{l|}{\cellcolor[HTML]{67FD9A}0.02345135} & \multicolumn{1}{l|}{\cellcolor[HTML]{67FD9A}0.02262303} & \cellcolor[HTML]{67FD9A}0.02239903 \\ \hline
\end{tabular} 
}
\caption{Mean Squared Error of the wN-W estimator across deterministic and random sampling designs, comparing Gaussian and Quartic kernels under different values of $\alpha$ and $H_{2}$ (with $H_{1} = 0.55$).}
\label{mse_tabla}
\end{table} 

\begin{figure}[H]
    \centering
    \begin{subfigure}[t]{0.48\textwidth}
        \centering
        \includegraphics[width=\textwidth]{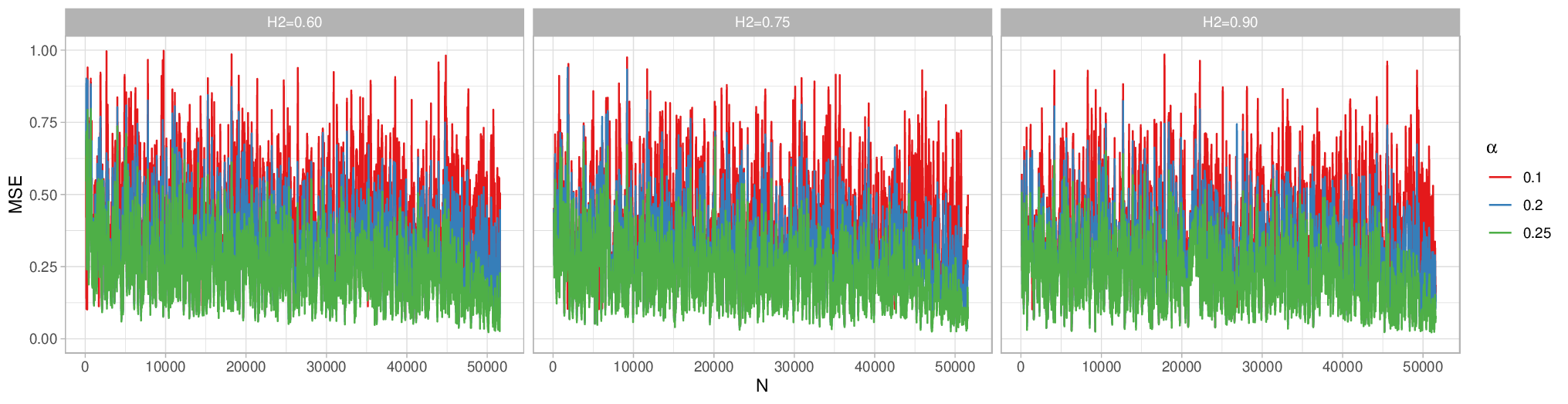}
        \caption{Deterministic - Gaussian}
    \end{subfigure}\hfill
    \begin{subfigure}[t]{0.48\textwidth}
        \centering
        \includegraphics[width=\textwidth]{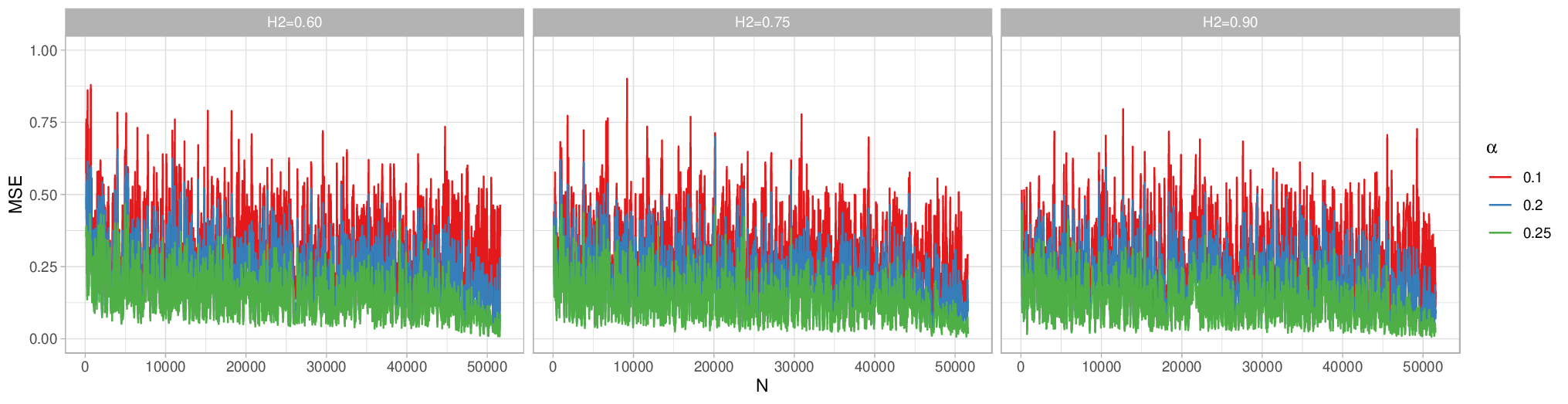}
        \caption{Deterministic - Quartic}
    \end{subfigure}
    
    \vspace{0.5cm} 
    
    \begin{subfigure}[t]{0.48\textwidth}
        \centering
        \includegraphics[width=\textwidth]{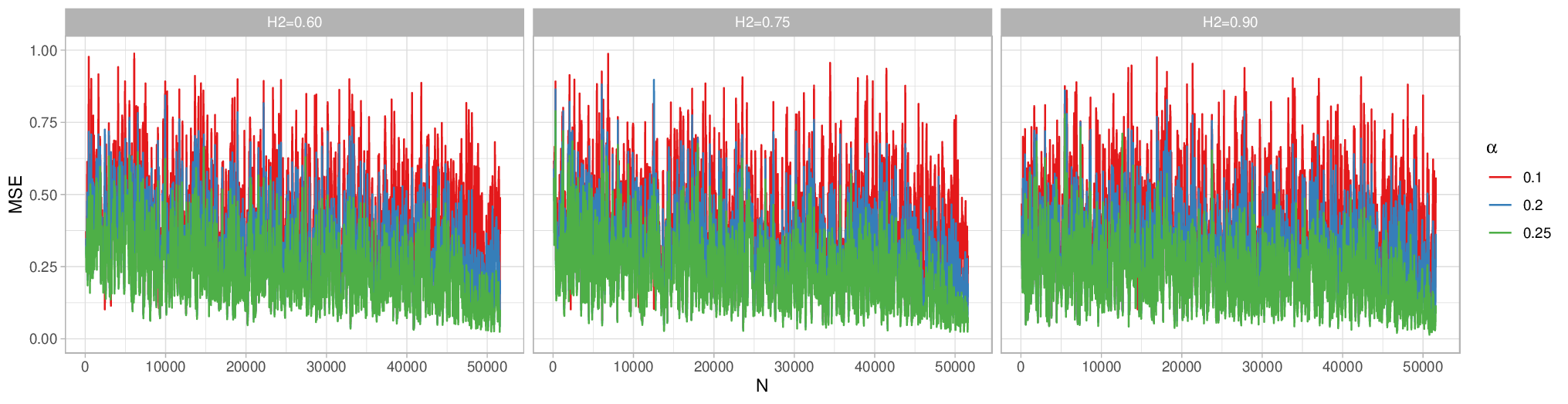}
        \caption{Random - Gaussian}
    \end{subfigure}\hfill
    \begin{subfigure}[t]{0.48\textwidth}
        \centering
        \includegraphics[width=\textwidth]{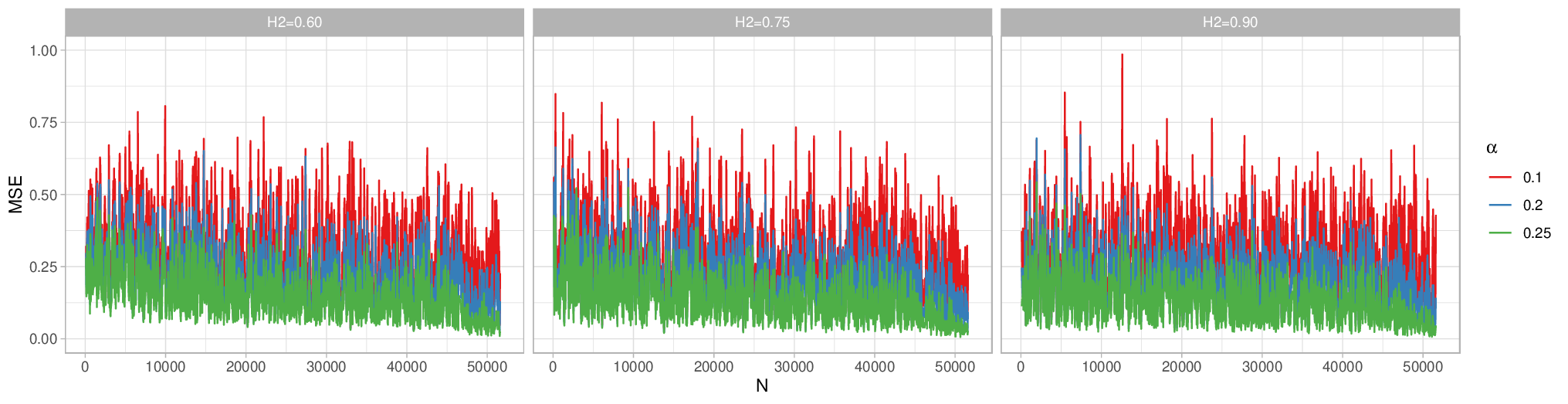}
        \caption{Random - Quartic}
    \end{subfigure}
\caption{MSE trajectories as a function of the expected number of observations ($N \in [100, 50000]$).} \label{mse-n}
\end{figure}

To further illustrate the asymptotic behavior of the MSE, Figure \ref{mse-n} displays the error trajectories for the evaluated methods as the expected number of observations $N$ increases from $100$ to $50.000$.
 
The empirical findings from Table \ref{mse_tabla} and Figure \ref{mse-n} reveal several key phenomena that closely align with our theoretical framework:
\begin{enumerate}
    \item Kernel Support Flexibility: Across all temporal designs and Hurst parameters, the bounded-support Quartic kernel consistently outperforms the Gaussian kernel, yielding MSE values approximately $1/2$ to $1/3$ lower. This empirical superiority reflects the less restrictive bandwidth constraints derived in Theorem \ref{thm:main} for bounded kernels.
    \item Influence of Memory and Bandwidth: The estimation error monotonically decreases as the noise persistence ($H_2$) approaches $1$. Concurrently, increasing the bandwidth parameter $\alpha$ (up to its valid upper bound, e.g., $\alpha=0.25$) significantly reduces the MSE. 
    \item Asymptotic Decay: Figure \ref{mse-n} demonstrates a strict convergence of the MSE toward zero as $N \to \infty$. Notably, for large samples ($N > 40000$) paired with $\alpha = 0.25$, the error trajectories stabilize and heavily contract, empirically validating the global $L^2$-consistency proven for both deterministic and random sampling contexts.
\end{enumerate}


\begin{rem}\label{extension}
It is noteworthy that the obtained results can be readily extended to the random time designated as ``jittered'', which represents random times as a regular partition subject to a bounded, positive random variable, thereby rendering the random times independent of one another. 

Likewise, the results obtained can be extended to stochastic processes having a similar covariance structure, such as the Hermite process.
\end{rem}

\begin{rem} \label{fourth-1}
Although almost sure convergence is not possible to prove, under the assumptions made for the possible values of the parameters $a$ and $b$, one can further establish convergence to higher orders, such as $L^{4}$ convergence.
\end{rem}

\section{Empirical Application: Vostok Ice Core Data} \label{real-data}
To illustrate and empirically validate the performance of the proposed estimator, paleoclimatic records extracted from the Vostok ice core in Antarctica are analyzed. The structured database consists of $n = 283$ joint empirical observations, where the explanatory covariate ($X_{\tau_i}$) corresponds to the atmospheric concentration of carbon dioxide ($CO_2$, expressed in parts per million by volume \footnote{\url{https://www.ncei.noaa.gov/pub/data/paleo/icecore/antarctica/vostok/co2nat.txt}.}) and the response variable ($Y_{\tau_i}$) represents the local temperature variation ($\Delta T$, in °C \footnote{\url{https://www.ncei.noaa.gov/pub/data/paleo/icecore/antarctica/vostok/deutnat.txt}.}) relative to the present value, derived from deuterium isotope analysis. Both measurements are collected at random observation times $\tau_i$, which denote the chronological age of the geological sample. Due to the natural dynamics of ice compression and extraction, the time series exhibits an inherently irregular sampling design, generating asymmetric time gaps ($\Delta \tau_i$) throughout the record. To ensure consistency with the topology of the spatial and temporal convergence theorems, the original time domain is linearly transformed to the standardized interval $[0, 1]$. \\

To theoretically validate the application of the estimator under irregular sampling, the method of moments was utilized to determine the parameter values governing the distribution of the random time gaps $\Delta \tau_i$. Based on the $n = 283$ normalized observations, the moment-based estimation yielded $a \approx 1.9098$ for the second moment and $b \approx 3.5713$ for the fourth moment. These results successfully verify the regularity conditions required for $L^2$ consistency, as they strictly satisfy the imposed theoretical constraints of $0 < a < 2$ and $0 < b < 2a$ (since $3.5713 < 3.8196$). Consequently, the empirical dataset intrinsically supports the mathematical assumptions of the proposed non-parametric model. \\

\begin{figure}[htbp]
\centering
\includegraphics[width = \textwidth]{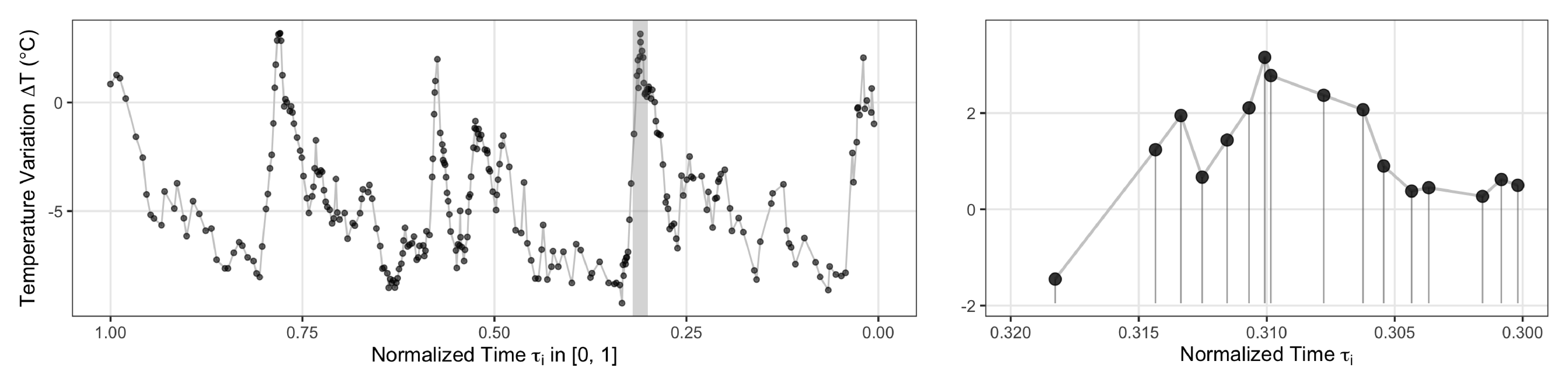}
\caption{Empirical observation times of the Vostok ice core dataset mapped to the normalized domain $\tau_i \in [0, 1]$.}
\label{serie}
\end{figure}

To motivate the necessity of the proposed non-parametric framework, Figure \ref{serie} displays the trajectory of the response variable over the standardized temporal domain $\tau_i \in [0, 1]$.

The left panel illustrates the complete macroscopic structure of the Vostok temperature variation record. The right panel provides a magnified view of a localized segment ($\tau_i \in [0.30, 0.32]$), explicitly exposing the intrinsic irregularity of the data acquisition process. The vertical projection lines demonstrate that the temporal distances between consecutive observations ($\Delta \tau_i$) are highly asymmetric, with periods of dense sampling immediately followed by substantial temporal gaps. This visual evidence underscores the inadequacy of classical equidistant estimation techniques and justifies the implementation of a mathematically rigorous temporal weighting scheme to achieve asymptotic consistency. \\

\begin{figure}[htbp]
\centering
\includegraphics[width = \textwidth]{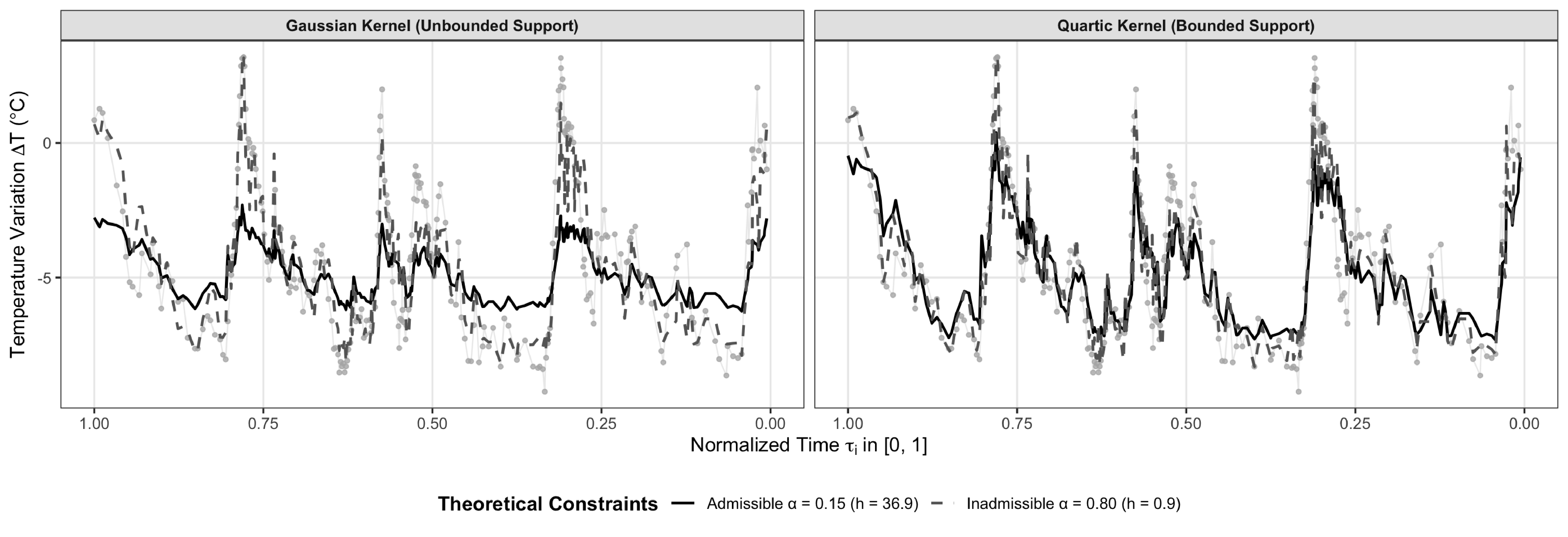}
\caption{Temporal reconstruction of the Vostok temperature variation using Bounded (Quartic) and Unbounded (Gaussian) kernels.}
\label{kernel-est}
\end{figure}

Figure \ref{kernel-est} illustrates the temporal reconstruction of the temperature variation ($\Delta T$) using the weighted Nadaraya-Watson estimator driven by the covariate $CO_2$. The estimation was executed over the standardized temporal domain $\tau_i \in [0, 1]$ utilizing both Bounded (Quartic) and Unbounded (Gaussian) kernels. 

To empirically evaluate the theoretical boundary conditions derived in the $L^2$ consistency theorems, two distinctive decay rates for the bandwidth sequence were contrasted: an admissible parameter ($\alpha = 0.15$) and an inadmissible parameter ($\alpha = 0.80$). Visually, the solid black line representing the admissible $\alpha$ demonstrates a robust capacity to capture the macroscopic structural trend of the glacial-interglacial cycles, effectively filtering the high-frequency stochastic noise inherent to the paleoclimatic record. Conversely, the dashed line corresponding to the inadmissible $\alpha$ exhibits severe structural instability and excessive variance, oscillating violently as it attempts to interpolate localized data points. From an asymptotic theory perspective, this behavior constitutes a classic manifestation of severe overfitting: by violating the upper bound of the regularity conditions, the inadmissible sequence forces the bandwidth $h$ to decay too rapidly, causing the estimator to memorize the random stochastic fluctuations rather than converging to the true underlying function. Therefore, the robust smoothing provided by the admissible parameter mathematically corroborates the necessity of the proposed moment constraints to guarantee the $L^2$ consistency and predictive generalization of the non-parametric model under irregular sampling. \\

In conclusion, the empirical application on the Vostok ice core dataset successfully bridges the gap between asymptotic theory and practical implementation for non-parametric regression under irregular sampling. The statistical evaluation of the temporal gaps ($\Delta \tau_i$) confirms that real-world paleoclimatic records inherently satisfy the strict moment conditions required for $L^2$ consistency. Furthermore, the graphical reconstruction underscores the critical importance of adhering to the theoretically derived bandwidth constraints; doing so avoids severe overfitting and ensures the accurate extraction of underlying macroscopic trends from highly noisy environments. Ultimately, the proposed weighted Nadaraya-Watson framework provides a mathematically rigorous and robust methodology for analyzing complex, non-equidistant time series, extending its potential applicability beyond paleoclimatology to broader contexts governed by stochastic observation times.

\section*{Acknowledgments}

Soledad Torres was partially supported by Fondecyt Reg. project number 1230807, 1221373, Mathamsud SMILE AMSUD230032, Mathamsud 
EXPLORE-SDE AMSUD240037 and Mathamsud SiJaVol AMSUD240024. 
This work was supported by Centro de Modelamiento Matemático (CMM) BASAL fund FB210005 for center of excellence from ANID-Chile.

\par

  \bibliographystyle{chicago}      
  \bibliography{no-p-bib}   

\end{document}